\documentclass[reqno]{amsart}
\usepackage{enumerate}
\usepackage[margin=1.2in]{geometry}
\usepackage{amssymb,url,color, booktabs,nccmath}
\definecolor{dgreen}{rgb}{0,0.5,0} 
\definecolor{revpurple}{rgb}{0.5,0,0.5} 
\usepackage{mathrsfs}
\usepackage{enumitem}
\usepackage{graphicx}

\usepackage{amsthm}
\usepackage{amsmath}

 \usepackage[colorlinks=true,linkcolor=blue!55!black,
  citecolor=blue!55!black,urlcolor=blue!65!black]{hyperref}
\usepackage{microtype}

\numberwithin{equation}{section}
\setlist{nosep,leftmargin=2em}

\usepackage{tikz}
\usetikzlibrary{arrows.meta,positioning,calc}

\newtheorem{theorem}{Theorem}[section]
\newtheorem{proposition}[theorem]{Proposition}
\newtheorem{lemma}[theorem]{Lemma}

\theoremstyle{definition}
\newtheorem{definition}[theorem]{Definition}

\theoremstyle{remark}
\newtheorem{remark}[theorem]{Remark}

\newcommand{\R}{\mathbb R}
\newcommand{\mR}{\mathbb R}
\newcommand{\E}{\mathbb E}
\newcommand{\mE}{\mathbb E}

\newcommand{\cB}{\mathcal B}
\newcommand{\Pp}{\mathbb P}
\newcommand{\1}{\mathbf 1}
\newcommand{\e}{\mathrm e}

\newcommand{\cV}{\mathcal V}
\newcommand{\dd}{\mathrm d}
\newcommand{\Ssym}{\mathbb S}
\newcommand{\Pucci}{\mathcal{P}^{-}_{\lambda,\Lambda}}
\newcommand{\qv}[1]{\langle #1 \rangle}
\newcommand{\cL}{\mathscr L}

\newcommand{\normK}[1]{\lVert #1\rVert_{\mathrm K}}
\newcommand{\fint}{\mathop{\rlap{\raisebox{0.18ex}{\scriptsize--}}\!\int}}
\newcommand{\dist}{\operatorname {dist}}

\DeclareMathOperator{\Law}{Law}
\DeclareMathOperator{\tr}{tr}
\DeclareMathOperator{\osc}{osc}

\title{Harnack estimates for nondivergence kinetic equations}
\author{Zimo Hao and Xicheng Zhang}

\address{Zimo Hao: School of Mathematics and Statistics, Beijing Institute of Technology, Beijing 100081, China,
Email: zmhao@bit.edu.cn}

\address{Xicheng Zhang:
School of Mathematics and Statistics, Beijing Institute of Technology, Beijing 100081, China\\
Faculty of Computational Mathematics and Cybernetics, Shenzhen MSU-BIT University, 518172 Shenzhen, China\\
Email: xczhang.math@bit.edu.cn
 }

\date{\today}

\begin{document}

\begin{abstract}
We establish a weak Harnack inequality for nonnegative strong supersolutions
and a Harnack inequality for nonnegative strong solutions of
the nondivergence-form kinetic equation
\[
 \partial_tu+v\cdot\nabla_xu+a(t,x,v):\nabla_v^2u
 +b(t,x,v)\cdot\nabla_vu+c(t,x,v)u=0,
\]
where the coefficients $a$, $b$, and $c$ are merely Borel measurable,
the matrix $a$ is uniformly elliptic, and $b$ and $c$ are bounded.
As a consequence, we derive interior kinetic H\"older estimates
for strong solutions.
Our approach develops a probabilistic kinetic analogue
of the Krylov--Safonov method.
The main new ingredient is a quantitative version of Krylov's estimate
for kinetic It\^o processes with progressively measurable coefficients,
which in turn yields the global existence of weak solutions
to the associated kinetic SDE with bounded measurable coefficients.
 These results appear to be the first Harnack and H\"older regularity
 theory for kinetic equations with merely measurable coefficients.

\end{abstract}
\maketitle
\enlargethispage{3pt}

\tableofcontents

\section{Introduction}

\subsection[Background]{Background}

The De Giorgi--Nash--Moser theory for divergence-form equations and the
Krylov--Safonov theory for nondivergence-form equations are two cornerstones
of the regularity theory of second-order elliptic and parabolic equations
\cite{DeGiorgi57,Nash58,Moser61,KrylovSafonov79,KrylovSafonov81,LSU68}.
Both theories show that bounded measurable leading coefficients still impose
quantitative regularity, but the mechanisms are fundamentally different.  In
divergence form, Caccioppoli inequalities for truncations combine with
Sobolev embedding and iteration.  In nondivergence form, no comparable energy
inequality is available; the Krylov--Safonov method instead uses
Aleksandrov-type measure estimates, growth lemmas, and quantitative hitting
probabilities for uniformly nondegenerate diffusions \cite{Krylov}.  

\vspace{1mm}

In this paper, we study the corresponding regularity problem for
nondivergence-form kinetic operators
\[
 \partial_t+v\cdot\nabla_x+a(t,x,v):\nabla_v^2,
\]
where the diffusion matrix $a$ is bounded, Borel measurable, and uniformly
elliptic.
Here $\partial_t+v\cdot\nabla_x$ describes free transport, while the
second-order diffusion acts only in the velocity variables.
The same transport--diffusion structure appears in kinetic models
such as the Landau equation \cite{Landau37,GIMV19}.

\vspace{1mm}

The
commutator identity
\[
 [\partial_{v_i},v\cdot\nabla_x]=\partial_{x_i}
\]
explains how velocity diffusion can nevertheless regularize the position
variable.  Thus, when the coefficients are smooth, the kinetic operator is
hypoelliptic by H\"ormander's theorem \cite{Hormander67}.
For the constant-coefficient, Kolmogorov had already
obtained an explicit fundamental solution in \cite{Kolmogorov34} by introducing the process $(\int_0^t W_s {\rm d}s, W_t)$, where $W$ is a standard $d$-dimensional Brownian motion.
The law of this process is invariant under the kinetic dilations
\[
  (t,x,v)\longmapsto(r^2t,r^3x,rv),
\]
which assign weights $2$, $3$, and $1$ to time, position, and
velocity, respectively.
Consistent with this scaling, Bouchut \cite{Bouchut2002}
established $L^p$-regularity estimates for the kinetic operator; see Chen and Zhang \cite{ChenZhang18} for the nonlocal setting.
For suitably regular coefficients, related heat-kernel and Harnack theory were investigated and developed in \cite{Polidoro97,PascucciPolidoro06Harnack,
BramantiBrandoliniLanconelliUguzzoni2010}.
Moreover, the intrinsic scaling also provides a framework for Schauder
estimates. Based on
anisotropic Besov spaces, a method for establishing such an estimate was developed in \cite{HaoWuZhang20} (see also \cite{HaoRocknerZhang26} the references therein).
For formulations and estimates in kinetic H\"older spaces adapted
to the transport flow and intrinsic dilations, see
\cite{HaoZhangZhuZhu24,ImbertSilvestre21Schauder}.

\vspace{1mm}

A rough-coefficient theory is well developed when the velocity diffusion is
in divergence form.  In that case,
$\operatorname{div}_v(a\nabla_v u)$ supplies the Caccioppoli inequalities
needed for truncation methods. Golse, Imbert, Mouhot, and Vasseur  \cite{GIMV19} 
established a kinetic
De Giorgi--Nash--Moser theory, including H\"older regularity and a Harnack
inequality, for weak solutions with rough coefficients.
Quantitative and geometric refinements, including weak Harnack inequalities,
kinetic Poincar\'e inequalities, and extensions to more general Kolmogorov
geometries, were obtained in
\cite{GuerandMouhot22,GI,
AnceschiRebucci22,NiebelZacher25}.  Every one of these
measurable-coefficient arguments uses the divergence structure at a decisive
stage.

\vspace{1mm}

In nondivergence form, for uniformly parabolic equations,
the classical Krylov--Safonov theory establishes Harnack inequalities
and interior H\"older estimates for bounded measurable,
uniformly elliptic diffusion coefficients
\cite{KrylovSafonov79,KrylovSafonov81,Krylov}.
Extending this theory to kinetic equations encounters two basic
obstacles: First, rewriting $a:\nabla_v^2u$ in divergence form
would introduce derivatives of the merely measurable matrix $a$,
so the standard energy argument underlying the divergence-form
theory cannot be applied directly except $a$ does not depend on velocity variable (see \cite{Zhang25McKeanVlasov}).
Second, the diffusion is degenerate in the position variables,
so the classical uniformly parabolic Krylov--Safonov theory
does not apply directly. These difficulties are reflected in the additional assumptions
imposed in existing Harnack results for nondivergence kinetic
operators: Abedin and Tralli
\cite{AbedinTralli19} proved Harnack inequalities assuming either uniform
continuity of the diffusion coefficients or the Cordes--Landis condition
$\Lambda/\lambda<1+1/(2d)$, where $\lambda$ and $\Lambda$
are the lower and upper ellipticity bounds of $a$, and
Hasanov's earlier one-dimensional result \cite{Hasanov12} also requires a
Cordes-type condition. Moreover, in \cite{Silvestre2023OpenProblems} Silvestre explicitly posed the kinetic Krylov--Safonov problem, conjecturing an interior H\"older estimate
for bounded classical solutions with constant and H\"older exponent
depending only on $d$, $\lambda$, and $\Lambda$.

\vspace{1mm}

The present paper closes this gap and gives a positive answer
to Silvestre's conjecture within the class of strong solutions. To the best of our knowledge, it establishes the first weak and strong Harnack inequalities, together with interior kinetic H\"older estimates, for nondivergence kinetic equations in arbitrary dimension with arbitrary bounded Borel measurable diffusion matrices satisfying uniform ellipticity. The equations also allow bounded Borel measurable velocity drifts and zeroth-order coefficients. These results provide a kinetic counterpart of the classical Krylov--Safonov theory in the class of strong solutions
(see Theorem~\ref{main} below).

\subsection{Main results}

Fix $d\in\mathbb N$, $0<\lambda\leqslant\Lambda<\infty$, and
${0\leqslant\hbar<\infty}$.  We consider the following nondivergence-form kinetic operator
\begin{equation}\label{eq:intro-operator}
 \cL_{a,b,c}u
 :=\partial_tu+v\cdot\nabla_xu+a(t,x,v):\nabla_v^2u
 +b(t,x,v)\cdot\nabla_vu+c(t,x,v)u,
\end{equation}
where $a:\R^{1+2d}\to\Ssym^d$,
$b:\R^{1+2d}\to\R^d$, and $c:\R^{1+2d}\to\R$ are Borel measurable and, for every
$(t,x,v)\in\R^{1+2d}$, satisfy
\begin{equation}\label{eq9}
 \lambda I_d\preceq a(t,x,v)\preceq\Lambda I_d,
 \quad |b(t,x,v)|, |c(t,x,v)|\leqslant \hbar.
\end{equation}
For notational simplicity, we write
$$
\Theta:=(d,\lambda,\Lambda,\hbar).
$$
For an open set $\Omega\subset\R^{1+2d}$, a \emph{strong solution}
(respectively, a \emph{strong supersolution}) of operator
$\cL_{a,b,c}$ in $\Omega$ is a function $w\in C(\Omega)$ whose
distributional derivatives satisfy
\[
 \partial_tw,\ \nabla_xw,\ \nabla_vw,\ \nabla_v^2w
 \in \cap_{p>1}L^p_{\mathrm{loc}}(\Omega),
\]
and for which
\[
 \cL_{a,b,c}w=0
 \quad\text{(respectively, $\cL_{a,b,c}w\leqslant0$)}
 \quad\text{Lebesgue-a.e. in }\Omega.
\]
We use the reference cylinders
\begin{align*}
 Q_+&:=\{(t,x,v): -1\leqslant t\leqslant0,\
                    |x|\leqslant1,\ |v|\leqslant1\},\\
 Q_-&:=\{(t,x,v): -10\leqslant t\leqslant-9,\
                    |x|\leqslant1,\ |v|\leqslant1\},
\end{align*}
and
\begin{align*}
Q_0&:=\{(t,x,v): |t|<1,\ |x|<1,\ |v|<1\},\\ 
Q'_0&:=\{(t,x,v): |t|<\tfrac14,\ |x|<\tfrac18,\ |v|<\tfrac12\}.
\end{align*}

The following theorem is the main result of the paper.

\begin{theorem}\label{main}
Let $Q\Subset\mathring Q_+$.  Under \eqref{eq9}, there exist
a bounded open set $\Omega$ satisfying
$Q_0\cup Q_+\cup Q_-\Subset\Omega$, exponents
$p_0,\alpha\in(0,1)$ depending only on $\Theta$, and positive constants
\[
 C_0=C_0(\Theta),\quad
 C_1=C_1\bigl(\Theta,d_{\rm K}(Q,Q_+^c)\bigr),\quad
 C_2=C_2(\Theta),
\]
where $d_{\rm K}$ is the kinetic distance defined in \eqref{eq:metric} below,
such that the following assertions hold:
\begin{itemize}
\item[(i)](Weak Harnack inequality) For every nonnegative strong supersolution $w$ of $\cL_{a,b,c}$ in
$\Omega$,
\begin{equation}\label{eq:intro-weak-harnack}
 \left(\fint_{Q_+}w^{p_0}\right)^{1/p_0}
 \leqslant C_0\inf_{Q_-}w.
\end{equation}
\item[(ii)](Strong Harnack inequality) For every nonnegative strong solution $w$ of
$\cL_{a,b,c}$ in $\Omega$,
\begin{equation}\label{eq:intro-strong-harnack}
 \sup_Qw\leqslant C_1\inf_{Q_-}w.
\end{equation}
\item[(iii)] (Interior H\"older estimate) For every strong solution $w$
of $\cL_{a,b,c}$ in $\Omega$,
\begin{equation}\label{eq:intro-holder}
 \sup_{z_1\ne z_2\in Q'_0}
 \frac{|w(z_1)-w(z_2)|}
      {d_{\rm K}(z_1,z_2)^\alpha}
 \leqslant C_2\|w\|_{L^\infty(Q_0)}.
\end{equation}
\end{itemize}
\end{theorem}

The estimates for the zero-order-free operator are proved in
Theorems~\ref{thm:weak-harnack}, \ref{thm:strong-harnack}, and
\ref{thm:holder}, and Section~\ref{sec:proof-main} considers the bounded
zeroth-order coefficient.  
Under a kinetic rescaling by $0<r\leqslant 1$, the lower-order coefficients become $r\,b$ and $r^2c$;
accordingly, the rescaled estimates depend only on the structural parameters
and the fixed reference geometry.

\vspace{1mm}

Our proof begins with establishing global and stopped Krylov-type occupation
estimates for kinetic It\^o processes with bounded progressively
measurable coefficients and uniformly elliptic velocity diffusion (see Theorems~\ref{thm:occupation} and
\ref{thm:stopped-occupation}).
The main difficulty is the absence of direct diffusion in the position
variables. Through It\^o's formula, these estimates yield a
kinetic Aleksandrov--Bakelman--Pucci type maximum estimate, obtained
probabilistically without invoking the classical nondegenerate ABP estimate.
More precisely, we consider the kinetic It\^o process
\begin{equation}\label{eq:intro-sde}
 \dd\mathsf X_s=\mathsf V_s\,\dd s,
 \qquad
 \dd\mathsf V_s=\boldsymbol\beta_s\,\dd s
 +\sqrt{2\mathbf A_s}\,\dd W_s,
\end{equation}
where $\mathbf A_s$ and $\boldsymbol\beta_s$ are arbitrary progressively
measurable processes and satisfy
\[
 \lambda I_d\preceq\mathbf A_s\preceq\Lambda I_d,
 \qquad
 |\boldsymbol\beta_s|\leqslant\hbar.
\]
In particular, they may be path dependent, and the
argument does not use a Markov semigroup or a pointwise heat-kernel bound.  In
the drift-free case, set
\[
 \mathsf Y_s:=\mathsf X_s-s\mathsf V_s,
 \qquad \mathsf Z_s:=(\mathsf Y_s,\mathsf V_s).
\]
We prove that, for Lebesgue-a.e. $s>0$, the law
$\mu_s:=\Law(\mathsf Z_s)$ admits a density $\rho_s$. Moreover, for some
$\digamma=\digamma(d,\lambda,\Lambda)\in(1,1+\frac{1}{2d})$ and every $T>0$,
we have
\begin{equation}\label{eq:intro-density-estimate}
 \int_0^T\|\rho_s\|_{L^\digamma(\R^{2d})}^{\digamma}\,\dd s
 \leqslant C_{d,\lambda,\Lambda}
 T^{1+2d-2d\digamma}.
\end{equation}
The main technical innovation of this part is a Gaussian
regularization that preserves the Fokker--Planck structure (see \eqref{eq:rg}):  it
closes the density estimate
under merely measurable coefficients and provides the foundation for the
Krylov estimates, weak existence, and Harnack inequalities.

\vspace{1mm}

As a consequence of \eqref{eq:intro-density-estimate}, 
we prove the existence of weak solutions for the kinetic SDE
$$
\dd\mathsf X_t=\mathsf V_t\dd t,
 \quad \dd\mathsf V_t=b(t,\mathsf X_t,\mathsf V_t)\dd t
 +\sqrt{2a(t,\mathsf X_t,\mathsf V_t)}\dd W_t,
$$
starting from every initial point $(\mathsf X_0,\mathsf V_0)=(x_0,v_0)\in\R^{2d}$; 
see Theorem~\ref{SDE}.  The Krylov estimates allow passage to the
limit with merely measurable coefficients via mollification.
The well-posedness results for regular coefficients refer to \cite{Zhang18, ChaudruDeRaynalMenozzi22}. Recently, in 
\cite{Zhang25McKeanVlasov}, the second named author established weak existence for density- and
distribution-dependent kinetic SDEs and well-posedness of generalized
martingale problems under a velocity-independent measurable diffusion
structure.

\vspace{1mm}

Next, we combine the Krylov-type occupation time
estimates with the kinetic geometry to show the Harnack inequalities. For this aim,
we introduce a \emph{link triple} $(K,B,O)$ (see Definition \ref{def:directed-link-triple} below), consisting of a source
set $K$, a target set $B$, and an outer domain $O$.
From every point of $K$, there is a controlled kinetic trajectory
that reaches the interior of $B$ at a later time, possibly followed
by a residence inside $B$.
The velocity controls have uniformly bounded energy, and all these
trajectories stay a uniform positive distance from the outer boundary.
This construction is adapted to the kinetic Lie-group structure
and anisotropic dilations (see \cite{LanconelliPolidoro94}).

\vspace{1mm}

The martingale tube estimate (Lemma~\ref{lem:tube}) then gives a
coefficient-uniform positive lower bound for the probability of
hitting $B$ before exiting $O$, uniformly over kinetic translations
and scales $0<r\leqslant1$ (see Lemma \ref{lem:uniform-link-hitting}).
When a positive residence interval is prescribed, the stopped Krylov
occupation estimate ensures that a uniform positive hitting bound
persists when $B$ is replaced by a measurable subset occupying a
sufficiently large fraction of it (see Lemma \ref{lem:nfh}).

\vspace{1mm}

The optional stopping theorem then converts these hitting estimates
into kinetic critical-density and double-ball estimates for nonnegative
strong supersolutions in Lemma \ref{lem:double-ball}. Their classical counterparts belong to the parabolic growth-theorem framework in \cite{FerrettiSafonov01}; an abstract
formulation of the critical-density and double-ball properties is given
in \cite{DFGL08}.

\vspace{1mm}

Combining these estimates with the kinetic ink-spots lemma
\cite[Corollary~10.2]{ImbertSilvestre20} (see also \cite{SchwabSilvestre2016}) yields the weak Harnack inequality for the case $c=0$ (Theorem \ref{thm:weak-harnack}). A further application of the Krylov estimate gives a local $L^{p}$--$L^\infty$ estimate for nonnegative strong solutions. Together with the weak Harnack inequality, this
establishes the strong Harnack inequality, Theorem \ref{thm:strong-harnack}, from which interior kinetic
H\"older estimates follow by oscillation decay (Theorem \ref{thm:holder}). Finally, a time-exponential transformation
handles the zeroth-order term $c$ in the weak Harnack estimate, while a
one-dimensional kinetic lift brings it into the strong Harnack and
H\"older estimates.

\subsection{Organization of the paper}

Section~2 proves the Krylov-type occupation estimate {and derives from it the}
near-full hitting estimate for link triples.  Section~3 constructs weak
solutions of the kinetic SDE with Borel coefficients and
provides the optional stopping theorem from the probabilistic estimates to
the PDE.  Section~4 proves the weak and strong Harnack inequalities for the
zero-order-free operator and derives the oscillation decay and interior
H\"older estimate.  Section~5 removes the zero-order restriction and proves
Theorem~\ref{main}.  Appendix~A collects the auxiliary stochastic estimates
and the kinetic ink-spots theorem with leakage.

\subsection{Notation}

We write $z=(t,x,v)$ and $z'=(s,y,w)$ for points of $\R^{1+2d}$.  The
kinetic group law and anisotropic dilations are defined by
\begin{equation}\label{eq:group}
 z\circ z'
 :=(t+s,x+y+sv,v+w),
 \qquad
 \delta_rz:=(r^2t,r^3x,rv),\quad r>0.
\end{equation}
The group inverse is
\[
 z^{-1}=(-t,-x+tv,-v).
\]
The homogeneous kinetic gauge is
\[
 \normK{z}:=\max\bigl\{|t|^{1/2},|x|^{1/3},|v|\bigr\},
\]
and its symmetric left-invariant distance is
\begin{equation}\label{eq:metric}
 d_{\rm K}(z,z')
 :=\max\left\{
   \normK{z^{-1}\circ z'},
   \normK{(z')^{-1}\circ z}
 \right\}.
\end{equation}
For $z_1,z_2,z_3\in\R^{1+2d}$, we have
\begin{equation}\label{eq:link-quasi-triangle}
 d_{\rm K}(z_1,z_3)
 \leqslant d_{\rm K}(z_1,z_2)+d_{\rm K}(z_2,z_3),
\end{equation}
and
\begin{equation}\label{Af3}
 \normK{z_1\circ z_2}
 \leqslant\normK{z_1}+\normK{z_2}.
\end{equation}
Thus $d_{\rm K}$ is a metric on $\R^{1+2d}$.  For nonempty sets
$E,F\subset\R^{1+2d}$, we use
\[
 d_{\rm K}(z,E):=\inf_{z'\in E}d_{\rm K}(z,z'),
 \quad
 d_{\rm K}(E,F):=\inf_{z\in E, z'\in F} d_{\rm K}(z,z').
\]
For $z=(t,x,v)$, $z'=(t',x',v')$, and $r>0$, define
\[
 \mathcal G_{z,r}(z'):=z\circ\delta_rz'
 =\left(t+r^2t',\,
        x+r^3x'+r^2t'v,\,
        v+rv'\right).
\]
Its inverse is
\[
 \mathcal G_{z,r}^{-1}(z')
 =
 \left(\frac{t'-t}{r^2},\,
       \frac{x'-x-(t'-t)v}{r^3},\,
       \frac{v'-v}{r}\right).
\]
For every Borel set $E\subset\R^{1+2d}$,
\begin{equation}\label{KiDi0}
 |\mathcal G_{z,r}(E)|=r^{4d+2}|E|.
\end{equation}
The maps $\mathcal G_{z,r}$ satisfy
\begin{equation}\label{KiDi2}
 \mathcal G_{z,r}\bigl(\mathcal G_{z',r'}(z'')\bigr)
 =\mathcal G_{\mathcal G_{z,r}(z'),rr'}(z''),
\end{equation}
and
\begin{equation}\label{KiDi1}
 d_{\rm K}\bigl(\mathcal G_{z,r}(z_1),
                 \mathcal G_{z,r}(z_2)\bigr)
 =r\,d_{\rm K}(z_1,z_2).
\end{equation}

We write $\Ssym^d$ for the space of real symmetric $d\times d$ matrices.
For $A,B\in\Ssym^d$, the notation $A\succeq B$ (equivalently,
$B\preceq A$) means that $A-B$ is positive semidefinite, while
$A\succ B$ means that $A-B$ is positive definite.  The Frobenius
contraction is defined by
\[
 A:B:=\tr(A^\top B)=\sum_{i,j=1}^dA_{ij}B_{ij},
\]
and the Frobenius norm is
\[
 |A|_{\rm F}:=\sqrt{A:A}.
\]
For $\Sigma\in\Ssym^d$ with $\Sigma\succ0$, let $g_\Sigma$ denote the density
of normal distribution $\mathcal N(0,2\Sigma)$:
\[
 g_\Sigma(x)
 :=\frac{1}{(4\pi)^{d/2}(\det\Sigma)^{1/2}}
 \exp\left(-\frac14x^\top\Sigma^{-1}x\right),
 \quad x\in\R^d.
\]

For a measurable set $E$, $|E|$ denotes its Lebesgue measure,
$\1_E$ its indicator, and
\[
 \fint_E f:=\frac1{|E|}\int_E f
\]
whenever $0<|E|<\infty$.  We use $E^c$ for the complement,
$\mathring E$ for the interior, and $E\Subset F$ to mean that
$\overline E$ is a compact subset of $\mathring F$.  Unless stated otherwise,
$C$ and $c$ denote positive constants whose values may change from line to
line; their parameter dependence is displayed when relevant.  We use
$\Theta=(d,\lambda,\Lambda,\hbar)$ for the structural parameter tuple.

\section{Probability estimates for general kinetic It\^o processes}
\label{sec:entropy}
{In this section, we derive probability estimates for general kinetic
It\^o processes that do not rely on the structure of any particular equation.}

Fix a filtered probability space $(\Omega,\mathcal F,(\mathcal F_t)_{t{\geqslant}0},\Pp)$
satisfying the usual conditions, and let $W$ be a $d$-dimensional
$(\mathcal F_t)$-Brownian motion.  
{Recall that $\Ssym^d$ denotes the space of real symmetric $d\times d$ matrices.}
Let
$\boldsymbol\beta:\R_+\times\Omega\to\R^d$ and
$\mathbf A:\R_+\times\Omega\to\Ssym^d$ be progressively measurable and assume
that, for some $\hbar<\infty$ and $0<\lambda<\Lambda<\infty$,
\begin{align}\label{eq:ellipticity}
 |\boldsymbol\beta_t|\leqslant \hbar,\qquad
 \lambda I_d\preceq \mathbf A_t\preceq\Lambda I_d.
\end{align}
{For brevity, write $\Theta:=(d,\lambda,\Lambda,\hbar)$.}
For fixed $x_0,v_0\in\R^d$, we define
\[
 \mathsf V_t:=v_0+\int_0^t\boldsymbol\beta_r\dd r
  +\int_0^t\sqrt{2\mathbf A_r}\dd W_r,
 \qquad
 \mathsf X_t:=x_0+\int_0^t\mathsf V_r\dd r.
\]
Thus $(\mathsf X_0,\mathsf V_0)=(x_0,v_0)$ and
\begin{equation}
 \dd\mathsf X_t=\mathsf V_t\dd t,
 \qquad
 \dd\mathsf V_t=\boldsymbol\beta_t\dd t+\sqrt{2\mathbf A_t}\dd W_t.
 \label{eq:controlled-sde}
\end{equation}
We refer to $(\mathsf X_t,\mathsf V_t)$ as a {\bf kinetic It\^o process}.

The aim of this section is to derive quantitative $L^p$-
estimates for the density {of the law of}
$(\mathsf X_t,\mathsf V_t)$.  These estimates may be viewed as a kinetic
analogue of the classical Krylov estimate \cite{Krylov}.
\subsection{A Krylov-type estimate}
{In this subsection we establish the crucial Krylov-type estimate for kinetic It\^o processes.}
We first consider the case
$x_0=v_0=0$ and $\boldsymbol\beta\equiv0$, and 
define $\mathsf Y_t:=\mathsf X_t-t\mathsf V_t$ and
\[
 \mathsf Z_t:=(\mathsf Y_t,\mathsf V_t)=
 \left(-\int_0^t s\sqrt{2\mathbf A_s}\dd W_s, \int_0^t\sqrt{2\mathbf A_s}\dd W_s\right),\quad
 \mu_t:=\Law(\mathsf Z_t).
\]
{For $B\in\cB(\R^{2d})$, define} the matrix-valued measure
\begin{equation}
 q_t(B):=\E\bigl[\mathbf A_t\1_{\{\mathsf Z_t\in B\}}\bigr].
 \label{eq:qmeasure}
\end{equation}
Then $q_t\ll\mu_t$.  For $t>0$, let
\begin{align}\label{eq:conditional-coefficient}
\mathfrak a_t(z):=\frac{\dd q_t}{\dd\mu_t}(z)
=\E\bigl[\mathbf A_t\mid\mathsf Z_t=z\bigr].
\end{align}
By \eqref{eq:ellipticity}, one has
\begin{equation}
 \lambda I_d\preceq\mathfrak a_t(z)\preceq\Lambda I_d
 \quad\text{for }\mu_t\text{-almost every }z.
 \label{eq:conditional-ellipticity}
\end{equation}
Define
\[
\mathscr D_{t,i}:=\partial_{v_i}-t\partial_{y_i},\qquad i=1,\dots,d,
\] 
and
\[
\mathscr D_t:=(\mathscr D_{t,1},\dots,\mathscr D_{t,d}).
\]
Since $\dd\mathsf Y_t=-t\,\dd\mathsf V_t$,  {It\^o's formula gives, for
$\varphi\in C_b^2(\R^{2d})$,}
\[
 \frac{\dd}{\dd t}\E\varphi(\mathsf Y_t,\mathsf V_t)
 =\E\left[(\mathbf A_t)_{ij}
  (\partial_{v_i}-t\partial_{y_i})
 (\partial_{v_j}-t\partial_{y_j})\varphi
  (\mathsf Y_t,\mathsf V_t)\right]
  {\quad\text{for Lebesgue-a.e. }t>0.}
\]
Here and below, repeated indices are summed.
\begin{lemma}
The curve $(\mu_t)_{t>0}$ satisfies, in the sense of distributions,
\begin{equation}
 \partial_t\mu_t={\mathscr D_{t}^2:q_t=\mathscr D_{t}^2:(\mathfrak a_t\mu_t)},
 \label{eq:sheared-fpe}
\end{equation}
\end{lemma}
In the following we write
\[
 \mathbf L_t:=(-tI_d, I_d)^\top\in\R^{2d\times d},
\]
where $\top$ denotes {the} transpose.  Then
\begin{align}\label{DD5}
\mathscr D_t=\mathbf L_t^\top\nabla_{y,v}.
\end{align}
{We will regularize $\mu_t$ by convolution with the Gaussian density of covariance matrix}
\[
 \Sigma_{s,t}:=\int_s^t\mathbf L_r\mathbf L_r^\top\dd r=
 \begin{pmatrix}
  \frac{t^3-s^3}{3}I_d&-\frac{t^2-s^2}{2}I_d\\[1mm]
  -\frac{t^2-s^2}{2}I_d&(t-s)I_d
 \end{pmatrix},\quad s<t.
\]
We have the following simple observation.
\begin{lemma}
\label{lem:leverage}
For any $s<t$, the matrix $\Sigma_{s,t}$ is invertible and
\begin{equation}
 \mathbf L_s^\top \Sigma_{s,t}^{-1}\mathbf L_s=\frac{4}{t-s}I_d,
 \qquad
 \mathbf L_s\mathbf L_s^\top\preceq \frac{4}{t-s}\Sigma_{s,t}.
 \label{eq:leverage}
\end{equation}
Moreover, for $\kappa{\geqslant}1$ and every $s,r>0$,
\begin{equation}
 \kappa\Sigma_{s,s+r}\preceq
\Sigma_{s,s+4\kappa r}.
 \label{eq:gram-compare}
\end{equation}

\end{lemma}
\begin{proof}
Let $r=t-s$.  {By a change of variables,} we have
\begin{align}\label{As8}
\Sigma_{s,t}=\Sigma_{s,s+r}=\int_0^r\mathbf L_{s+u}\mathbf L_{s+u}^\top\dd u
=\mathcal S_s\left(\int_0^r\mathbf L_{u}\mathbf L_{u}^\top\dd u\right) \mathcal S^\top_s,
\end{align}
where
\[
\mathcal S_s=
 \begin{pmatrix}
  I_d&-sI_d\\[1mm]
  0&I_d
 \end{pmatrix},\quad
 \mathcal S_s^{-1}=
 \begin{pmatrix}
  I_d&sI_d\\[1mm]
 0&I_d
 \end{pmatrix}.
\]
Hence, it suffices to consider $s=0$.  In this case
\[
 \Sigma_{0,r}^{-1}
 =\begin{pmatrix}
  12r^{-3}I_d&6r^{-2}I_d\\
  6r^{-2}I_d&4r^{-1}I_d
 \end{pmatrix},
\]
which gives the first identity in \eqref{eq:leverage}; the second follows by
$$
\mathcal S_s^{-1}\mathbf L_s(S_s^{-1}\mathbf L_s)^\top=
\begin{pmatrix}
  0_d&0_d\\
  0_d&I_d
 \end{pmatrix}
 \preceq 4
 \begin{pmatrix}
  \frac{r^2}{3}I_d&-\frac{r}{2}I_d\\[1mm]
  -\frac{r}{2}I_d&I_d
 \end{pmatrix}=
 \frac{4}{r}\int_0^r\mathbf L_{u}\mathbf L_{u}^\top\dd u.
$$
Finally, by \eqref{As8}, \eqref{eq:gram-compare} reduces to the positivity of
\[
 \begin{pmatrix}
  ((4\kappa)^3-\kappa)r^3/3&-((4\kappa)^2-\kappa)r^2/2\\
  -((4\kappa)^2-\kappa)r^2/2&3\kappa r
 \end{pmatrix}.
\]
Its determinant is $\kappa^2(8\kappa-5/4)r^4>0$.
This completes the proof.
\end{proof}
For $0\leqslant s<t$ and $\sigma>0$, define
\begin{equation}
 \begin{aligned}
 p^\sigma_{s,t}
 :=\mu_s*g_{\sigma \Sigma_{s,t}},\quad
 q^\sigma_{s,t}
 :=q_s*g_{\sigma \Sigma_{s,t}},
 \end{aligned}
 \label{eq:regularized}
\end{equation}
where $g_{\sigma \Sigma_{s,t}}$ is the density of $\mathcal N(0,2\sigma \Sigma_{s,t})$, and
\[
 \overline{\mathfrak a}_{s,t}^{\sigma}
 :=q_{s,t}^\sigma/p_{s,t}^\sigma.
\]
By \eqref{eq:ellipticity}, we also have
\begin{equation}
 \lambda I_d\preceq\overline{\mathfrak a}_{s,t}^\sigma
 \preceq\Lambda I_d.
 \label{eq:posterior-elliptic}
\end{equation}
\begin{lemma}
Let $\digamma\in(1,2)$ and $\sigma\in(0,\lambda)$.
For every $t>0$, the map
$s\mapsto\int_{\R^{2d}}(p_{s,t}^\sigma)^\digamma\dd z$
is locally absolutely continuous on $[0,t)$.
Moreover, for Lebesgue-a.e. $s\in(0,t)$,
\begin{equation}
 \frac{\dd}{\dd s}\int_{\mathbb R^{2d}}(p_{s,t}^\sigma)^\digamma
 \leqslant \frac{d\digamma(\digamma-1)}{4(\lambda-\sigma)}
      \sum_{i,j} \int_{\R^{2d}} (p_{s,t}^\sigma)^\digamma|\mathscr D_{s,j}
       \overline{\mathfrak a}_{s,t}^{\sigma,ij}|^2\dd z.
 \label{eq:entropy-upper}
\end{equation}
\end{lemma}
\begin{proof}
Since  $\partial_s\Sigma_{s,t}=-\mathbf L_s\mathbf L_s^\top$
{and $\mathscr D_s=\mathbf L_s^\top\nabla_{y,v}$,} we have
\[
 \partial_sg_{\sigma \Sigma_{s,t}}
 =-\sigma{\mathscr D_{s}^2}: (g_{\sigma \Sigma_{s,t}}I_d),
\]
which, together with \eqref{eq:sheared-fpe}, implies that
\begin{align}
\partial_s p^\sigma_{s,t}
 =\partial_s\mu_s*g_{\sigma \Sigma_{s,t}}+\mu_s*\partial_sg_{\sigma \Sigma_{s,t}}
 ={\mathscr D_s^2}:
 \left((\overline{\mathfrak a}_{s,t}^{\sigma}
              -\sigma I_d)p_{s,t}^\sigma\right).\label{eq:rg}
\end{align}
Define
\[
 \mathfrak w_i:=\mathscr D_{s,i}\log p_{s,t}^\sigma,
 \qquad
 \mathfrak g_i:=\mathscr D_{s,j}
       \overline{\mathfrak a}_{s,t}^{\sigma,ij}.
\]
By integration by parts, we have
\begin{equation}
 \frac{\dd}{\dd s}\int(p_{s,t}^\sigma)^\digamma
 =-\digamma(\digamma-1)\int (p_{s,t}^\sigma)^\digamma
 \left[\mathfrak w^\top(\overline{\mathfrak a}_{s,t}^\sigma-\sigma I_d)
 \mathfrak w+\mathfrak g\cdot\mathfrak w\right].
 \label{eq:entropy-exact}
\end{equation}
By \eqref{eq:posterior-elliptic}, Young's inequality, and
$|\sum_j b_j|^2\leqslant d\sum_j|b_j|^2$, we further have
\begin{align*}
 \frac{\dd}{\dd s}\int_{\R^{2d}}(p_{s,t}^\sigma)^\digamma\dd z
 &\leqslant \frac{\digamma(\digamma-1)}{4(\lambda-\sigma)}
       \int_{\R^{2d}}(p_{s,t}^\sigma)^\digamma|\mathfrak g|^2\dd z
       \leqslant \frac{d\digamma(\digamma-1)}{4(\lambda-\sigma)}
      \sum_{i,j}\int_{\R^{2d}}(p_{s,t}^\sigma)^\digamma
      |\mathscr D_{s,j}\overline{\mathfrak a}_{s,t}^{\sigma,ij}|^2\dd z.
\end{align*}
This completes the proof.
\end{proof}
For $\Sigma\succ0$ and $s>0$, we define
\[
 \cV_s^\digamma(\Sigma):=\sum_{i,j=1}^d\int_{\R^{2d}}
 \frac{|q_s^{ij}*g_\Sigma|^2}
      {(\mu_s*g_\Sigma)^{2-\digamma}}\dd z,
 \qquad 1<\digamma<2.
\]
\begin{lemma}
\label{lem:data-processing}
If $\Sigma_2\succeq\Sigma_1\succ0$, then
\[
 \cV_s^\digamma(\Sigma_2)\leqslant\cV_s^\digamma(\Sigma_1).
\]
Moreover, for every $\Sigma\succ0$,
\begin{equation}
 0\leqslant\cV_s^\digamma(\Sigma)
 \leqslant d\Lambda^2\int_{\R^{2d}}(\mu_s*g_\Sigma)^\digamma\dd z.
 \label{eq:variance-budget}
\end{equation}
\end{lemma}

\begin{proof}
Suppose that $\Sigma_2-\Sigma_1\succ0$ and write
$\Sigma_2=\Sigma_1+\Sigma$.  By the Gaussian semigroup property,
\[
 q_s^{ij}*g_{\Sigma_2}
 =q_s^{ij}*g_{\Sigma_1}*g_\Sigma=:u_s^{ij}*g_\Sigma,
 \qquad
 \mu_s*g_{\Sigma_2}
 =\mu_s*g_{\Sigma_1}*g_\Sigma=:p_s*g_\Sigma.
\]
By H\"older's inequality, we have
\[
|u_s^{ij}*g_\Sigma|^2
 \leqslant
((|u_s^{ij}|^2/p_s^{2-\digamma})*g_\Sigma)\,(p_s^{2-\digamma}*g_\Sigma)
\]
and
$$
|q_s^{ij}*g_\Sigma|^2\leqslant
(|\mathfrak a_s^{ij}|^2*g_\Sigma)\,(\mu_s*g_\Sigma)^2.
$$
Hence,
$$
 \frac{|u_s^{ij}*g_\Sigma|^2}{(p_s*g_\Sigma)^{2-\digamma}}
 \leqslant
 \left(\frac{|u_s^{ij}|^2}{p_s^{2-\digamma}}\right)*g_\Sigma,\quad
\frac{|q_s^{ij}*g_\Sigma|^2}{(\mu_s*g_\Sigma)^{2-\digamma}}\leqslant
(|\mathfrak a_s^{ij}|^2*g_\Sigma)\,(\mu_s*g_\Sigma)^\digamma.
$$
Integrating and summing over $i,j$, we obtain the monotonicity and \eqref{eq:variance-budget}.
\end{proof}

\begin{lemma}[de Bruijn-type inequality]
\label{lem:debruijn}
For fixed $0<s<t$ and $\digamma\in(1,2)$, we have
\begin{equation}
\begin{split}
 -\partial_\sigma\cV_s^\digamma(\sigma\Sigma_{s,t})
 \geqslant \frac{2(2-\digamma)}{\digamma}
 \sum_{i,j}\int_{\R^{2d}}(p_{s,t}^\sigma)^\digamma
 (\nabla\overline{\mathfrak a}_{s,t}^{\sigma,ij})^\top
 \Sigma_{s,t}\nabla\overline{\mathfrak a}_{s,t}^{\sigma,ij}\dd z.
\end{split}
 \label{eq:debruijn}
\end{equation}
\end{lemma}
\begin{proof}
{Fix $s<t$, write $\Sigma:=\Sigma_{s,t}$, and set }
\[
 p:=\mu_s*g_{\sigma\Sigma},\qquad
 q^{ij}:=q_s^{ij}*g_{\sigma\Sigma}.
\]
Then
\[
 \partial_\sigma p=\Sigma:\nabla^2p,\qquad
 \partial_\sigma q^{ij}=\Sigma:\nabla^2q^{ij}.
\]
Set $\overline{\mathfrak a}^{ij}:=q^{ij}/p$.  Differentiating and integrating by parts gives
\begin{align}\label{Aw3}
\begin{split}
 \partial_\sigma\int_{\R^{2d}}\frac{|q^{ij}|^2}{p^{2-\digamma}}\dd z
 ={}&-2\int_{\R^{2d}}p^\digamma
        (\nabla \overline{\mathfrak a}^{ij})^\top\Sigma\nabla \overline{\mathfrak a}^{ij}\dd z\\
 &-4(\digamma-1)\int_{\R^{2d}}p^\digamma \overline{\mathfrak a}^{ij}
        (\nabla \overline{\mathfrak a}^{ij})^\top\Sigma\nabla\log p\dd z\\
 &-\digamma(\digamma-1)\int_{\R^{2d}}p^\digamma|\overline{\mathfrak a}^{ij}|^2
        (\nabla\log p)^\top\Sigma\nabla\log p\dd z.
\end{split}
\end{align}
Note that
\begin{align*}
 &2(\nabla \overline{\mathfrak a}^{ij})^\top\Sigma\nabla \overline{\mathfrak a}^{ij}
 +4(\digamma-1)\overline{\mathfrak a}^{ij}(\nabla \overline{\mathfrak a}^{ij})^\top\Sigma\nabla\log p +\digamma(\digamma-1)|\overline{\mathfrak a}^{ij}|^2
       (\nabla\log p)^\top\Sigma\nabla\log p\\
 &=\frac{2(2-\digamma)}{\digamma}
       (\nabla \overline{\mathfrak a}^{ij})^\top\Sigma\nabla \overline{\mathfrak a}^{ij}
 +\digamma(\digamma-1)
 \left(\overline{\mathfrak a}^{ij}\nabla\log p+\frac{2}{\digamma}\nabla \overline{\mathfrak a}^{ij}\right)^\top
 \Sigma
 \left(\overline{\mathfrak a}^{ij}\nabla\log p+\frac{2}{\digamma}\nabla \overline{\mathfrak a}^{ij}\right).
\end{align*}
Hence,
\begin{align*}
 \partial_\sigma\int_{\R^{2d}}\frac{|q^{ij}|^2}{p^{2-\digamma}}\dd z \leqslant{}&-\frac{2(2-\digamma)}{\digamma}
       \int_{\R^{2d}}p^\digamma
        (\nabla \overline{\mathfrak a}^{ij})^\top\Sigma\nabla \overline{\mathfrak a}^{ij}\dd z.
\end{align*}
{Summing over $i,j$ gives \eqref{eq:debruijn}. }
\end{proof}
\begin{remark}
The classical de Bruijn equality states that the derivative
of Shannon entropy along Gaussian convolution is the Fisher
information (see \cite{Stam1959}). More precisely, for $f\geqslant0$,
\[
\frac{\dd}{\dd\sigma}\mathcal H(f*g_{\sigma\Sigma})
=
\int_{\R^{2d}}(f*g_{\sigma\Sigma})
\left|\Sigma^{1/2}\nabla\log (f*g_{\sigma\Sigma})\right|^2\,\dd z,
\]
where
$\mathcal H(f):=-\int_{\R^{2d}}f\log f\,\dd z$.
The above \eqref{Aw3} is a variant of the classical de Bruijn equality
(see \cite[Proposition 1]{Chafai2004}).
\end{remark}

We are now ready to prove the following key density estimate.

\begin{theorem}\label{thm:entropy}
For Lebesgue-a.e. $t>0$, the measure $\mu_t$ is absolutely continuous with
respect to Lebesgue measure.  Moreover, there exists $\digamma=\digamma(d,\lambda,\Lambda)\in(1,1+\frac{1}{2d})$
such that, for every $T>0$,
\begin{align}\label{Dual1}
 \int_0^T\|\dd\mu_t/\dd z\|_{L^\digamma}^\digamma\dd t
 \leqslant C_{d,\lambda,\Lambda}\,T^{1-2d(\digamma-1)}.
\end{align}
\end{theorem}
\begin{proof}
Note that
\begin{align*}
 \sum_{i,j}|\mathscr D_{s,j}
       \overline{\mathfrak a}_{s,t}^{\sigma,ij}|^2
 \stackrel{\eqref{DD5}}{\leqslant}\sum_{i,j}(\nabla\overline{\mathfrak a}_{s,t}^{\sigma,ij})^\top
       \mathbf L_s\mathbf L_s^\top
       \nabla\overline{\mathfrak a}_{s,t}^{\sigma,ij}
\stackrel{\eqref{eq:leverage}}{\leqslant}\frac4{t-s}
       \sum_{i,j}(\nabla\overline{\mathfrak a}_{s,t}^{\sigma,ij})^\top
       \Sigma_{s,t}\nabla\overline{\mathfrak a}_{s,t}^{\sigma,ij}.
\end{align*}
By \eqref{eq:entropy-upper} and \eqref{eq:debruijn}, for
$\sigma\in(0,\lambda)$ and Lebesgue-a.e. $s\in(0,t)$,
\[
 \frac{\dd}{\dd s}\int_{\R^{2d}}(p_{s,t}^\sigma)^\digamma\dd z
 \leqslant
 -\frac{d\digamma^2(\digamma-1)}
 {2(2-\digamma)(\lambda-\sigma)(t-s)}
 \partial_\sigma\cV_s^\digamma(\sigma\Sigma_{s,t}).
\]
Integrating both sides in $\sigma$ from $\sigma_0:=\lambda/4$ to $\sigma_1:=\lambda/2$, we have
\begin{align*}
 &\frac{\dd}{\dd s}\int_{\sigma_0}^{\sigma_1}\!\!\!
       \int_{\R^{2d}}(p_{s,t}^\sigma)^\digamma\dd z\,\dd\sigma\leqslant
 \frac{d\digamma^2(\digamma-1)}{(2-\digamma)\lambda}
 \frac{\cV_s^\digamma(\sigma_0\Sigma_{s,t})
       -\cV_s^\digamma(\sigma_1\Sigma_{s,t})}{t-s}.
\end{align*}
Furthermore, for $\varepsilon\in(0,t)$,  integrating {in} $s$ from $0$ to $t-\varepsilon$, we get
\begin{align*}
 &\int_{\sigma_0}^{\sigma_1}\!\!\!
       \int_{\R^{2d}}(p_{t-\varepsilon,t}^\sigma)^\digamma\dd z\,\dd\sigma
       -\int_{\sigma_0}^{\sigma_1}\!\!\!
       \int_{\R^{2d}}(p_{0,t}^\sigma)^\digamma\dd z\,\dd\sigma
       \leqslant
 \frac{d\digamma^2(\digamma-1)}{(2-\digamma)\lambda}
 \int^{t-\varepsilon}_0\frac{\cV_s^\digamma(\sigma_0\Sigma_{s,t})
       -\cV_s^\digamma(\sigma_1\Sigma_{s,t})}{t-s}\dd s.
\end{align*}
Finally, for $\varepsilon\in(0,T)$, 
integrating in $t$ from
$\varepsilon$ to $T$ and applying Fubini's theorem together with the change of variables $r=t-s$, 
we obtain
\begin{equation}\label{eq:terminal-entropy}
\begin{aligned}
 &\int_\varepsilon^T\!\!\!\int_{\sigma_0}^{\sigma_1}\!\!\!
       \int_{\R^{2d}}|p_{t-\varepsilon,t}^\sigma|^\digamma
       \dd z\,\dd\sigma\,\dd t
       -\int_\varepsilon^T\!\!\!\int_{\sigma_0}^{\sigma_1}\!\!\!
       \int_{\R^{2d}}|p_{0,t}^\sigma|^\digamma
       \dd z\,\dd\sigma\,\dd t\\
 &\quad\leqslant
 \frac{d\digamma^2(\digamma-1)}{(2-\digamma)\lambda}
 \int_0^{T-\varepsilon}\!\!\!\int_\varepsilon^{T-s}
 \frac{\cV_s^\digamma(\sigma_0\Sigma_{s,s+r})
       -\cV_s^\digamma(\sigma_1\Sigma_{s,s+r})}{r}
       \dd r\,\dd s.
\end{aligned}
\end{equation}
For notational simplicity, we write
\[
 F^{\sigma_0}_s(r):=\cV_s^\digamma(\sigma_0\Sigma_{s,s+r}),\qquad r>0.
\]
Since $\sigma_1/\sigma_0=2$, by \eqref{eq:gram-compare} we have
\[
 \sigma_0\Sigma_{s,s+r}\preceq\sigma_1\Sigma_{s,s+r}
 \preceq\sigma_0\Sigma_{s,s+8r}.
\]
Thus, by Lemma \ref{lem:data-processing},
\[
F^{\sigma_0}_s(r)= \cV_s^\digamma(\sigma_0\Sigma_{s,s+r})
 \geqslant\cV_s^\digamma(\sigma_1\Sigma_{s,s+r})
 \geqslant\cV_s^\digamma(\sigma_0\Sigma_{s,s+8r})=F^{\sigma_0}_s(8r)\geqslant 0.
\]
Consequently, for $R>\varepsilon$,
\begin{equation}\label{eq:log-telescope}
\begin{aligned}
 0&\leqslant\int_\varepsilon^R
 \frac{\cV_s^\digamma(\sigma_0\Sigma_{s,s+r})
       -\cV_s^\digamma(\sigma_1\Sigma_{s,s+r})}{r}\dd r
       \leqslant\int_\varepsilon^R\frac{F^{\sigma_0}_s(r)-F^{\sigma_0}_s(8r)}r\dd r\\
&=\int_\varepsilon^{8\varepsilon}\frac{F^{\sigma_0}_s(r)}r\dd r
   -\int_R^{8R}\frac{F^{\sigma_0}_s(r)}r\dd r\leqslant F^{\sigma_0}_s(\varepsilon)
   \int_\varepsilon^{8\varepsilon}\frac{\dd r}r
\leqslant F^{\sigma_0}_s(\varepsilon)\log 8.
\end{aligned}
\end{equation}
Hence,
\begin{align}\label{Ag1}
\int_0^{T-\varepsilon}\!\!\!\int_\varepsilon^{T-s}
 \frac{\cV_s^\digamma(\sigma_0\Sigma_{s,s+r})
       -\cV_s^\digamma(\sigma_1\Sigma_{s,s+r})}{r}
       \dd r\,\dd s\leqslant\log 8\int_0^{T-\varepsilon}F^{\sigma_0}_s(\varepsilon)\dd s.
\end{align}
By \eqref{eq:variance-budget}, for $\sigma\in[\sigma_0,\sigma_1]$, we have
\[
 F^{\sigma_0}_s(\varepsilon)
 \leqslant d\Lambda^2\int_{\R^{2d}}
             |p_{s,s+\varepsilon}^{\sigma_0}|^\digamma\dd z
             \leqslant d\Lambda^2 2^{d\digamma}\int_{\R^{2d}}
             |p_{s,s+\varepsilon}^{\sigma}|^\digamma\dd z,
\]
where we have used the Gaussian density estimate
\[
 g_{\sigma_0\Sigma}
 \leqslant\left(\frac{\sigma}{\sigma_0}\right)^d g_{\sigma\Sigma}
 \leqslant2^d g_{\sigma\Sigma}.
\]
In particular,
\begin{align}
 \int_0^{T-\varepsilon}F^{\sigma_0}_s(\varepsilon)\dd s
 &\leqslant
 \frac{d\Lambda^2\,2^{d\digamma}}{\sigma_1-\sigma_0}
 \int_0^{T-\varepsilon}\int_{\sigma_0}^{\sigma_1}
       \int_{\R^{2d}}|p_{s,s+\varepsilon}^\sigma|^\digamma
       \dd z\,\dd\sigma\,\dd s\nonumber\\
 &=
 \frac{4d\Lambda^2\,2^{d\digamma}}{\lambda}
 \int_\varepsilon^T\int_{\sigma_0}^{\sigma_1}
       \int_{\R^{2d}}|p_{t-\varepsilon,t}^\sigma|^\digamma
       \dd z\,\dd\sigma\,\dd t.\label{Ag2}
\end{align}
On the other hand, since
\[
 \|p_{0,t}^\sigma\|_\infty\leqslant C_{d,\lambda}t^{-2d},
 \quad \int_{\R^{2d}}p_{0,t}^\sigma\dd z=1,
\]
we have
\begin{equation}\label{eq:initial-entropy}
 \int_{\R^{2d}}|p_{0,t}^\sigma|^\digamma\dd z
 \leqslant\|p_{0,t}^\sigma\|_\infty^{\digamma-1}
 \leqslant C_{d,\lambda,\digamma}t^{-2d(\digamma-1)}.
\end{equation}
We now choose $\digamma\in(1,2)$
sufficiently close to $1$ so that
\[
2d(\digamma-1)<1,\quad \frac{4d^2\Lambda^2\digamma^2(\digamma-1)\,2^{d\digamma}\log8}
 {(2-\digamma)\lambda^2}
 \leqslant\frac12.
\]
Combining \eqref{eq:terminal-entropy}, \eqref{Ag1}, \eqref{Ag2} and  \eqref{eq:initial-entropy}, we obtain
\begin{align*}
&\int_\varepsilon^T\!\!\!\int_{\sigma_0}^{\sigma_1}\!\!\!
       \int_{\R^{2d}}|p_{t-\varepsilon,t}^\sigma|^\digamma
       \dd z\,\dd\sigma\,\dd t\leqslant
 \int_\varepsilon^T\!\!\!\int_{\sigma_0}^{\sigma_1}\!\!\!
       \int_{\R^{2d}}|p_{0,t}^\sigma|^\digamma
       \dd z\,\dd\sigma\,\dd t\\
&\quad       +\frac{4d^2\Lambda^2\digamma^2(\digamma-1)\,2^{d\digamma}\log8}
 {(2-\digamma)\lambda^2}
       \int_\varepsilon^T\!\!\!\int_{\sigma_0}^{\sigma_1}\!\!\!
       \int_{\R^{2d}}|p_{t-\varepsilon,t}^\sigma|^\digamma
       \dd z\,\dd\sigma\,\dd t\\
 &\qquad\leqslant C_{d,\lambda,\digamma}\int^T_\varepsilon t^{-2d(\digamma-1)}\dd t+
 \frac12\int_\varepsilon^T\!\!\!\int_{\sigma_0}^{\sigma_1}\!\!\!
       \int_{\R^{2d}}|p_{t-\varepsilon,t}^\sigma|^\digamma
       \dd z\,\dd\sigma\,\dd t.
\end{align*}
Hence, for all $\varepsilon\in(0,T)$,
\[
 \int_\varepsilon^T\!\!\!\int_{\sigma_0}^{\sigma_1}\!\!\!
 \int_{\R^{2d}}(p_{t-\varepsilon,t}^\sigma)^\digamma
 \dd z\,\dd\sigma\,\dd t
 \leqslant \frac{2C_{d,\lambda,\digamma}T^{1-2d(\digamma-1)}}{1-2d(\digamma-1)}.
\]
Let $f\in C_c([0,T]\times\R^{2d})$ and $\xi\sim\mathcal N(0, I_{2d})$. 
For each $\sigma$ and $t>0$, we have
\begin{align*}
\mE f(t,\mathsf Z_t)=\lim_{\varepsilon\to 0}\mE f\left(t,\mathsf Z_t+\sigma\sqrt{\Sigma_{t-\varepsilon,t}} \xi\right)
=\lim_{\varepsilon\to 0}\int_{\R^{2d}} f(t,z) p_{t-\varepsilon,t}^\sigma(z)\dd z
\end{align*}
Hence, by Fatou's lemma and H\"older's inequality,
\begin{align*}
\int^T_0|\mE f(t,\mathsf Z_t)|\dd t
&\leqslant\varliminf_{\varepsilon\to 0}\frac{1}{\sigma_1-\sigma_0} \int_\varepsilon^T\!\!\!\int_{\sigma_0}^{\sigma_1}\!\!\!
\int_{\R^{2d}} f(t,z) p_{t-\varepsilon,t}^\sigma(z)\dd z\dd t\\
&\leqslant\tfrac{4}{\lambda}\|f\|_{L^{\digamma_*}([0,T]\times\R^{2d})}(C_{d,\lambda,\Lambda}T^{1-2d(\digamma-1)})^{1/\digamma},
\end{align*}
where $\digamma_*=\frac{\digamma}{\digamma-1}$.
It follows that, for Lebesgue-a.e. $t>0$, the law $\mu_t$ of $\mathsf Z_t$ admits a density
and that \eqref{Dual1} holds.
\end{proof}

We can now prove the following general Krylov-type estimate.
\begin{theorem}[Krylov's estimate]\label{thm:occupation}
Under \eqref{eq:ellipticity}, for any $T>0$, there {exist} $p=p(d,\lambda,\Lambda)>1+2d$ and  
{a} constant $C_0=C_0(\Theta, T)>0$
 such that{,} for every  $f\in L^p([0,T]\times\R^{2d})$,
\begin{equation}
 \E\left(\int_0^T f(t,\mathsf X_t,\mathsf V_t)\dd t\right)
 \leqslant \min\left\{T\|f\|_{L^\infty([0,T]\times\R^{2d})},C_0\|f\|_{L^p([0,T]\times\R^{2d})}
 \right\}.
 \label{eq:unit-occupation}
\end{equation}
\end{theorem}
\begin{proof}
Since $\E\left(\int_0^T f(t,\mathsf X_t,\mathsf V_t)\dd t\right)
 \leqslant T\|f\|_{L^\infty([0,T]\times\R^{2d})}$ is trivial, we only need to show
 $$
 \E\left(\int_0^T f(t,\mathsf X_t,\mathsf V_t)\dd t\right)
 \leqslant C_0\|f\|_{L^p([0,T]\times\R^{2d})}.
 $$
Let $\digamma\in(1,1+\frac1{2d})$ be as in Theorem \ref{thm:entropy},
and let $\digamma^*:=\frac{\digamma}{\digamma-1}$.
For $\boldsymbol\beta\equiv0$, {it follows from the change of variables, Theorem \ref{thm:entropy},
and H\"older's inequality that}
\begin{align}
&\E\left(\int_0^T f(t,\mathsf X_t,\mathsf V_t)\dd t\right)
=\int_0^T\!\!\!\int_{\R^{2d}} f(t,x_0+y+t(v_0+v),v_0+v)\mu_t(\dd y,\dd v)\dd t\nonumber\\
&\quad\leqslant\left(\int^T_0\|f(t,x_0+\cdot+t(v_0+\cdot),v_0+\cdot)\|_{\digamma^*}^{\digamma^*}\dd t\right)^{1/\digamma^*}
\left(\int^T_0\|\dd\mu_t/\dd z\|^\digamma_{L^\digamma}\dd t\right)^{1/\digamma}\nonumber\\
&\quad\leqslant \|f\|_{L^{\digamma^*}([0,T]\times\R^{2d})}(C_{d,\lambda,\Lambda}\,T^{1+2d-2d\digamma})^{1/\digamma}.
\label{DD8}
\end{align}
For general $\boldsymbol\beta$, {we use Girsanov's theorem}. 
Fix $T>0$.  Write $\sigma_t:=\sqrt{2\mathbf A_t}$ and
$\theta_t:=\sigma_t^{-1}\boldsymbol\beta_t$.  Since
$|\theta_t|^2\leqslant \hbar^2/(2\lambda)$, by Novikov's condition, we may
define a probability measure $\mathbb Q$ by
\[
 \frac{\dd\mathbb Q}{\dd\Pp}
 =\exp\left(
 -\int_0^T\theta_t\cdot\dd W_t
 -\frac12\int_0^T|\theta_t|^2\dd t
 \right).
\]
Under $\mathbb Q$,
\[
 W_t^{\mathbb Q}:=W_t+\int_0^t\theta_s\dd s,\qquad 0\leqslant t\leqslant T,
\]
is a Brownian motion, and therefore
\[
 \dd\mathsf X_t=\mathsf V_t\dd t,\qquad
 \dd\mathsf V_t=\sqrt{2\mathbf A_t}\dd W_t^{\mathbb Q}.
\]
For $q>1$ and $q^*:=q/(q-1)$, by H\"older's inequality, we have
\begin{align}
 \E\left(\int_0^T f(t,\mathsf X_t,\mathsf V_t)\dd t\right)
{\leqslant}
 \left(\E_{\mathbb Q}
 \left(\frac{\dd\Pp}{\dd\mathbb Q}\right)^q\right)^{1/q}
 \left(\E_{\mathbb Q}\left|\int_0^T f(t,\mathsf X_t,\mathsf V_t)\dd t\right|^{q_*}\right)^{1/q_*}.\label{DD9}
\end{align}
By the exponential martingale and the bound
$|\theta_t|^2\leqslant \hbar^2/(2\lambda)$, we have
\begin{equation}\label{DD7}
 \left(\E_{\mathbb Q}
 \left(\frac{\dd\Pp}{\dd\mathbb Q}\right)^q\right)^{1/q}
 \leqslant
 \exp\left(\frac{(q-1)\hbar^2T}{4\lambda}\right).
\end{equation}
Combining \eqref{DD8}, \eqref{DD9} and \eqref{DD7}, we obtain the desired estimate with $p=q^*\digamma^*$.
\end{proof}

\begin{remark}
{It remains an open question whether \eqref{eq:unit-occupation} holds for all $p>1+2d$.}
\end{remark}
\subsection{Stopped occupation time estimates}
In this subsection, we use Theorem \ref{thm:occupation} to derive a
stopped occupation-time estimate needed later in the proof of the Harnack inequality.
{For fixed $z_0:=(t_0,x_0,v_0)\in{\R}^{1+2d}$, define} the augmented process
$$
\mathbf Z_s^{\rm aug}{:=}(t_0+s, \mathsf X_s,\mathsf V_s),\quad s\geqslant 0,\ \ \mathbf Z_0^{\rm aug}=z_0.
$$
{We write $\Pp_{z_0}$ and $\E_{z_0}$ for probability and expectation,
respectively, for the augmented process started at $z_0$.}
{For a Borel set $D\subset\R^{1+2d}$ and an open set
$O\subset\R^{1+2d}$, define the hitting and exit times by}
\[
 {T_D:=\inf\{s\geqslant0:\mathbf Z_s^{\rm aug}\in D\},\qquad}
 \tau_O:=\inf\{s{\geqslant}0:\mathbf Z_s^{\rm aug}\notin O\}.
\]
{We use the convention that the infimum of the empty set is $+\infty$.}
{We have the following stopped estimate.}
\begin{theorem} \label{thm:stopped-occupation}
Let $K\subset\R^{1+2d}$ be a nonempty bounded Borel set, and let
$O\subset\R^{1+2d}$ be a bounded open set such that $\overline K\Subset O$.
Define
\begin{equation}
 H_0:=\max\left\{1,
  \sup_{(t,x,v)\in O}t
  -\inf_{(t,x,v)\in K}t\right\}.
 \label{eq:fixed-shape-time-span}
\end{equation}
Under \eqref{eq:ellipticity}, there {exist} $p=p(d,\lambda,\Lambda)>1+2d$ and  
{a} constant $C_0=C_0(\Theta, H_0)>0$ such that 
for all $(z,r)\in\R^{1+2d}\times(0,1]$ and for every Borel set $E\subset\R^{1+2d}$,
\begin{equation}
 \sup_{z_0\in \mathcal G_{z,r}(K)}
 \E_{z_0}\left(\int_0^{\tau_{\mathcal G_{z,r}(O)}}
 \1_E(\mathbf Z_t^{\rm aug})\dd t\right)
 \leqslant C_0 r^2\left(
 \frac{|E\cap \mathcal G_{z,r}(O)|}{r^{4d+2}}
 \right)^{\frac1p}.
 \label{eq:stopped-occupation}
\end{equation}
\end{theorem}

\begin{proof}
We prove the estimate in two steps.  

\medskip
\noindent\textit{Step 1: a unit-scale stopped estimate.}
Let $E\subset\R^{1+2d}$ be an arbitrary Borel set.
Fix $z_0=(t_0,x_0,v_0)\in K$. {Let $\tau_O$ be the exit time of the
augmented process $\mathbf Z^{\rm aug}$ from $O$ when started at $z_0$.}
By definition of $\tau_O$, we have
\[
 \tau_O\leqslant \sup_{(t,x,v)\in O}t-t_0\leqslant H_0,
\]
and for any $s\in[0,\tau_O)$,
$$
\1_E(\mathbf Z_s^{\rm aug})
=\1_{E\cap O}(\mathbf Z_s^{\rm aug}).
$$ 
Hence, by \eqref{eq:unit-occupation} with $f(t,z):=\1_{E\cap O}(z)$,
\[
\mE_{z_0}\left(\int_0^{\tau_O}\1_E(\mathbf Z_s^{\rm aug})\dd s\right)
 \leqslant
\mE_{z_0}\left( \int_0^{H_0}\1_{E\cap O}(\mathbf Z_s^{\rm aug})\dd s\right)
\leqslant C_0|E\cap O|^{1/p}.
\]

\medskip
\noindent\textit{Step 2: scaling the process.}
Now fix $z=(t_z,x_z,v_z)$ and $r\in{(0,R_0]}$.
For $z_0=(t_0,x_0,v_0)\in \mathcal G_{z,r}(K)$, set
\[
 \bar z_0:=\mathcal  G_{z,r}^{-1}(t_0,x_0,v_0)
 =
 \left(
 \frac{t_0-t_z}{r^2},\,
 \frac{x_0-x_z-(t_0-t_z)v_z}{r^3},\,
 \frac{v_0-v_z}{r}
 \right)\in K.
\]
For $s\geqslant0$, define
\[
 \widehat{\mathbf Z}_s^{\rm aug}
 :=\mathcal  G_{z,r}^{-1}(\mathbf Z_{r^2s}^{\rm aug})
 =\bigl(\widehat t_s,{\widehat{\mathsf X}_s},{\widehat{\mathsf V}_s}\bigr),
\]
where
\[
 \widehat t_s=\frac{t_0-t_z}{r^2}+s,\qquad
 {\widehat{\mathsf X}_s}
 =
 \frac{{\mathsf X_{r^2s}}-x_z-(t_0+r^2s-t_z)v_z}{r^3},
 \qquad
 {\widehat{\mathsf V}_s}=\frac{{\mathsf V_{r^2s}}-v_z}{r}.
\]
In particular, if we define
$$
\widehat W_s:=r^{-1}W_{r^2s},\quad
{\widehat{\mathbf A}_s:=\mathbf A_{r^2s}},\quad
 {\widehat{\boldsymbol\beta}_s:=r\boldsymbol\beta_{r^2s}},
$$
{then the rescaled process satisfies}
\[
 \dd{\widehat{\mathsf X}_s}={\widehat{\mathsf V}_s}\dd s,\quad
 \dd{\widehat{\mathsf V}_s}={\widehat{\boldsymbol\beta}_s}\dd s
 +\sqrt{2{\widehat{\mathbf A}_s}}\dd\widehat W_s.
\]
Note that
\begin{align}\label{As2}
\lambda I_d\preceq{\widehat{\mathbf A}_s}\preceq\Lambda I_d,\quad
|{\widehat{\boldsymbol\beta}_s}|\leqslant r\hbar\leqslant \hbar.
\end{align}
Let
\[
 \widehat\tau_O
 :=\inf\{s{\geqslant}0:\widehat{\mathbf Z}_s^{\rm aug}\notin O\}.
\]
Then
\[
 \tau_{\mathcal G_{z,r}(O)}=\inf\left\{s{\geqslant}0:\mathbf Z_s^{\rm aug}\notin \mathcal G_{z,r}(O)\right\}
 =\inf\left\{s{\geqslant}0:\widehat{\mathbf Z}_{s/r^2}^{\rm aug}\notin O\right\}=r^2\widehat\tau_O.
\]
Set
\[
 \widehat E:=\mathcal  G_{z,r}^{-1}(E\cap \mathcal G_{z,r}(O))=\mathcal  G_{z,r}^{-1}(E)\cap O\subset O.
\]
{Using the change of variables $s=r^2t$, the bounds \eqref{As2}, and
the estimate from Step 1, we obtain}
\begin{align*}
 \E_{z_0} \int_0^{\tau_{\mathcal G_{z,r}(O)}}
 \1_E(\mathbf Z_s^{\rm aug})\dd s
 &=
 r^2\E_{\bar z_0}\int_0^{\widehat\tau_O}
 \1_{\widehat E}(\widehat{\mathbf Z}_t^{\rm aug})\dd t\leqslant C_0
|\widehat E|^{1/p}.
\end{align*}
{Since
$|\widehat E|=r^{-4d-2}|E\cap\mathcal G_{z,r}(O)|$, this is exactly
\eqref{eq:stopped-occupation}.}
\end{proof}

\subsection{{Hitting probabilities for measurable sets}}

{In this subsection, we derive estimates for the probability that a
kinetic It\^o process hits a measurable set before exiting a domain.  We first
introduce the notion of a link triple of subsets of $\R^{1+2d}$.}
Given a starting point $z=(t_z,x_z,v_z)$ and an absolutely continuous
velocity curve $\phi:{\R}_+\to\R^d$ with $\phi(0)=0$, define its
deterministic lift to the augmented space by
\begin{equation}
 \boldsymbol\gamma^{z,\phi}_s
 :=\left(t_z+s,
 x_z+sv_z+\int_0^s\phi(u)\dd u,v_z+\phi(s)\right),
 \quad s\geqslant 0.
 \label{eq:deterministic-lift}
\end{equation}
{For $h>0$ and an open set $O$, define the distance of this path from
$O^c$ by}
\begin{equation}
d_{\rm K}\left(\boldsymbol\gamma^{z,\phi}_{[0,h]},O^c\right):=
\inf_{0\leqslant s\leqslant h}
 d_{\rm K}(\boldsymbol\gamma^{z,\phi}_s,O^c).
 \label{eq:path-clearance}
\end{equation}

{We now introduce the notion of a link triple.}

\begin{definition}[Link triple]\label{def:directed-link-triple}
Fix parameters $\eta_*, A_*, h_*\in(0,\infty)$ and $\kappa_*\geqslant 0$.
{Let $K,B,O\subset\R^{1+2d}$ be such that $K$ is nonempty and compact,
$B$ is compact with nonempty interior, and $O$ is bounded and open, with}
\[
 K\cup B\Subset O.
\]
We call $(K,B,O)$ a link triple if
for every $z=(t_z,x_z,v_z)\in K$, there exist
$h\in{(0,h_*]}$ and an absolutely continuous centered velocity control
$\phi:[0,h+\kappa_*]\to\R^d$ satisfying $\phi(0)=0$ and
\[
\int_0^{h+\kappa_*}|\dot\phi(s)|^2\dd s\leqslant A_*,\quad d_{\rm K}\left(\boldsymbol\gamma^{z,\phi}_{[0,h+\kappa_*]},O^c\right)\geqslant2\eta_*,\quad
d_{\rm K}\left(\boldsymbol\gamma^{z,\phi}_{[h,h+\kappa_*]},B^c\right)\geqslant2\eta_*.
\]
The class of all such link triples is denoted by
$\mathcal A_{\rm link}(\eta_*,A_*,h_*,\kappa_*)$.
\end{definition}
{The following figure depicts the link triple constructed in
Lemma~\ref{lem:elementary-directed-steering}; horizontal placement represents
the time coordinate only schematically.}
\begin{center}
\begin{tikzpicture}[
  x=0.94cm,y=0.94cm,
  >=Stealth,
  every node/.style={font=\scriptsize,text=black},
  steering/.style={very thick,blue!65!black,-{Stealth[length=5pt]}},
  residence/.style={very thick,red!70!black,-{Stealth[length=5pt]}},
  timeguide/.style={densely dashed,gray!70}
]
  \node[font=\small\bfseries] at (6.5,4.68)
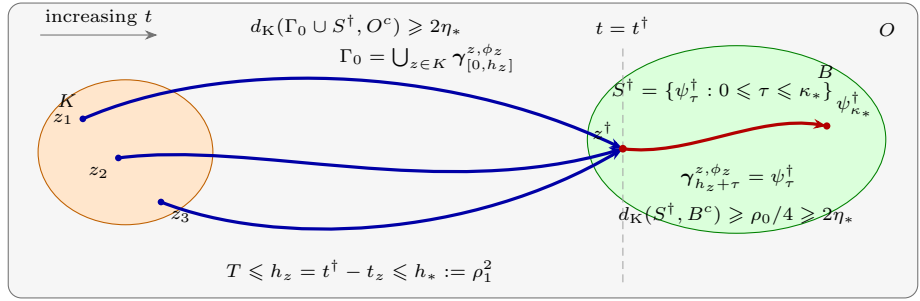
    {Link triple in Lemma~\ref{lem:elementary-directed-steering}};
  \draw[fill=gray!7,draw=black!55,rounded corners=5pt]
    (0.1,0.2) rectangle (12.9,4.35);
  \node[anchor=north east] at (12.72,4.20) {$O$};

  \draw[->,black!55] (0.55,3.90) -- (2.20,3.90)
    node[midway,above] {increasing $t$};
  \node at (5.0,4.03)
    {$d_{\rm K}(\Gamma_0\cup S^\dagger,O^c)\geqslant2\eta_*$};

  \draw[fill=orange!20,draw=orange!75!black]
    (1.75,2.25) ellipse (1.23 and 1.02);
  \node at (0.93,2.98) {$K$};
  \coordinate (zone) at (1.15,2.72);
  \coordinate (ztwo) at (1.65,2.17);
  \coordinate (zthree) at (2.25,1.55);
  \fill[blue!65!black] (zone) circle (1.25pt)
    node[left] {$z_1$};
  \fill[blue!65!black] (ztwo) circle (1.25pt)
    node[below left] {$z_2$};
  \fill[blue!65!black] (zthree) circle (1.25pt)
    node[below right] {$z_3$};

  \coordinate (target) at (8.75,2.30);
  \coordinate (resend) at (11.62,2.62);
  \draw[fill=green!14,draw=green!50!black]
    (10.35,2.43) ellipse (2.10 and 1.32);
  \node at (11.60,3.38) {$B$};

  \draw[steering] (zone)
    .. controls (3.05,3.55) and (6.20,3.52) .. (target);
  \draw[steering] (ztwo)
    .. controls (3.55,2.42) and (6.20,1.55) .. (target);
  \draw[steering] (zthree)
    .. controls (3.85,0.92) and (6.60,1.02) .. (target);
  \node[fill=gray!7,inner sep=1pt] at (6,3.58)
    {$\Gamma_0=\bigcup_{z\in K}
      \boldsymbol\gamma^{z,\phi_z}_{[0,h_z]}$};
  \node at (5.05,0.58)
    {$T\leqslant h_z=t^\dagger-t_z\leqslant h_*:=\rho_1^2$};

  \draw[timeguide] (8.75,0.42) -- (8.75,3.78);
  \node[anchor=south] at (8.75,3.78) {$t=t^\dagger$};
  \fill[red!70!black] (target) circle (1.45pt)
    node[above left] {$z^\dagger$};
  \draw[residence] (target)
    .. controls (9.65,2.18) and (10.65,2.78) .. (resend);
  \fill[red!70!black] (resend) circle (1.35pt)
    node[above right] {$\psi^\dagger_{\kappa_*}$};
  \node at (10.15,3.12)
    {$S^\dagger=\{\psi^\dagger_\tau:0\leqslant\tau\leqslant\kappa_*\}$};
  \node at (10.38,1.90)
    {$\boldsymbol\gamma^{z,\phi_z}_{h_z+\tau}=\psi^\dagger_\tau$};
  \node at (10.35,1.38)
    {$d_{\rm K}(S^\dagger,B^c)\geqslant\rho_0/4\geqslant2\eta_*$};

  \node[align=center] at (6.50,-0.22)
    {Every $z\in K$ is steered to the common point $z^\dagger$;
     thereafter all paths coincide with
     $S^\dagger\subset\operatorname{int}B$.};
\end{tikzpicture}
\end{center}

\begin{lemma}[Existence of link triples]
\label{lem:elementary-directed-steering}
Let $K\subset \R^{1+2d}$ be a nonempty compact set and
$B\subset \R^{1+2d}$ be compact with nonempty interior. Assume that
there is a point
$z^\dagger=(t^\dagger,x^\dagger,v^\dagger)\in B$ satisfying
\begin{align}\label{Af1}
 T:=t^\dagger-\sup_K t>0\quad\text{and}\quad
\rho_0:=d_{\rm K}(z^\dagger,B^c)>0,\quad \rho_1:=\sup_{z\in K,z'\in B}d_{\rm K}(z,z').
\end{align}
Then one can choose a bounded open set $O$ such that
$(K,B,O)\in\mathcal A_{\rm link}(\eta_*,A_*,h_*,\kappa_*)$ for some
$\eta_*,A_*,h_*,\kappa_*\in(0,\infty)$ depending only on $d$, $T$, and $\rho_0,\rho_1$.
\end{lemma}
\begin{proof}
For $z=(t_z,x_z,v_z)\in K$, set
\[
 h_z:=t^\dagger-t_z,\quad
 y_z:=x^\dagger-x_z-h_zv_z,\quad
 w_z:=v^\dagger-v_z.
\]
By \eqref{eq:metric}, we have
\begin{equation}\label{eq:link-intrinsic-bounds}
 T\leqslant h_z\leqslant\rho_1^2=:h_*,
 \qquad |y_z|\leqslant\rho_1^3,
 \qquad |w_z|\leqslant\rho_1.
\end{equation}
For $0\leqslant s\leqslant h_z$, define
\begin{equation}\label{eq:link-control}
 \phi_z(s)
 :=\left(\frac{6y_z}{h_z^2}-\frac{2w_z}{h_z}\right)s
 +\left(\frac{3w_z}{h_z^2}-\frac{6y_z}{h_z^3}\right)s^2.
\end{equation}
By direct calculation, we have
\[
 \phi_z(0)=0,\quad
 \phi_z(h_z)=w_z,\quad
 \int_0^{h_z}\phi_z(s)\dd s=y_z.
\]
Let ${\boldsymbol\gamma^{z,\phi_z}_s}$ be defined by \eqref{eq:deterministic-lift}.
Consequently,
\begin{equation}\label{eq:link-reaches-target}
 \boldsymbol\gamma^{z,\phi_z}_0=z,
 \qquad
 \boldsymbol\gamma^{z,\phi_z}_{h_z}=z^\dagger.
\end{equation}
Moreover, \eqref{eq:link-control} and
\eqref{eq:link-intrinsic-bounds} imply that
\begin{equation}\label{eq:link-action-bound}
 \int_0^{h_z}|\dot\phi_z(s)|^2\dd s
 \leqslant 16\left(\frac{|y_z|^2}{h_z^3}
              +\frac{|w_z|^2}{h_z}\right)
 \leqslant 16\left(\frac{\rho_1^6}{T^3}
              +\frac{\rho_1^2}{T}\right).
\end{equation}
We next construct the residence part of the path quantitatively.  Let
\[
 \kappa_*:=\left(\frac{\rho_0}{4}\right)^2,
 \qquad
 \psi^\dagger_\tau
 :=(t^\dagger+\tau,x^\dagger+\tau v^\dagger,v^\dagger),
 \quad 0\leqslant\tau\leqslant\kappa_*.
\]
Note that
$d_{\rm K}(z^\dagger,\psi^\dagger_\tau)=\sqrt\tau$.  For every
$z'\in B^c$, by \eqref{eq:link-quasi-triangle}, we have
\[
 \rho_0
 \leqslant d_{\rm K}(z^\dagger,z')
 {\leqslant \sqrt\tau+d_{\rm K}(\psi^\dagger_\tau,z')},
\]
and hence
\begin{equation}\label{eq:link-residence-clearance}
 {d_{\rm K}}(S^\dagger,B^c)
 \geqslant\tfrac{\rho_0}{4},
 \qquad
 S^\dagger:=\{\psi^\dagger_\tau:0\leqslant\tau\leqslant\kappa_*\}.
\end{equation}
In particular, 
$$
S^\dagger\subset\operatorname{int}B.
$$
Now, extend the control in \eqref{eq:link-control} by
\[
 \phi_z(s):=w_z,
 \qquad h_z\leqslant s\leqslant h_z+\kappa_*.
\]
The extension is absolutely continuous because $\phi_z(h_z)=w_z$, and
\begin{equation}\label{eq:link-residence-path}
 \boldsymbol\gamma^{z,\phi_z}_{h_z+\tau}=\psi^\dagger_\tau,
 \qquad 0\leqslant\tau\leqslant\kappa_*.
\end{equation}
It has zero additional action.  Thus, the required action is bounded by $A_*$ for every $z\in K${, where}
\[
 A_*:=1+16\left(\frac{\rho_1^6}{T^3}
                 +\frac{\rho_1^2}{T}\right).
\]
{Since $h_z$, $y_z$, and $w_z$ depend continuously on $z\in K$ and
$h_z\geqslant T$,} the set
\[
 \Gamma_0
 :=\{\boldsymbol\gamma^{z,\phi_z}_{\theta h_z}:
       z\in K,\ 0\leqslant\theta\leqslant1\}
\]
is compact, being the continuous image of $K\times[0,1]$.  Hence
$\Gamma_0\cup S^\dagger\cup B$ is compact.  Define
\[
 O:=\left\{z\in\R^{1+2d}:
       \inf_{z'\in \Gamma_0\cup S^\dagger\cup B}d_{\rm K}(z,z')<1\wedge(\tfrac{\rho_0}{4})\right\}.
\]
Then $O$ is bounded and
{open, and $K\cup B\Subset O$.}  Since every controlled path is contained
in $\Gamma_0\cup S^\dagger$, one also has
\begin{equation}\label{eq:link-outer-clearance}
{d_{\rm K}}\left(
 \boldsymbol\gamma^{z,\phi_z}_{[0,h_z+\kappa_*]},O^c
 \right)\geqslant1\wedge(\tfrac{\rho_0}{4}),
 \qquad z\in K.
\end{equation}
Finally, choose
\[
 \eta_*:=\tfrac12
 \left(1\wedge(\tfrac{\rho_0}{4})\right)>0.
\]
The bounds \eqref{eq:link-intrinsic-bounds},
\eqref{eq:link-action-bound}, \eqref{eq:link-residence-clearance}, and
\eqref{eq:link-outer-clearance}, together with
\eqref{eq:link-residence-path}, verify all the requirements of Definition
\ref{def:directed-link-triple}.  The displayed choices show that
$\eta_*,A_*,h_*$ and $\kappa_*$ depend only on
$d,T,\rho_0,$ and $\rho_1$.
\end{proof}

{We now prove uniform hitting and occupation-time estimates for link triples.}
\begin{lemma}[Uniform hitting and occupation time estimates]
\label{lem:uniform-link-hitting}
Let $(K,B,O)\in\mathcal A_{\rm link}(\eta_*,A_*, h_*,\kappa_*)$.
Under \eqref{eq:ellipticity}, there is a constant $\ell_0=\ell_0(\Theta,\eta_*,A_*,h_*,\kappa_*)>0$
such that for any $z\in \R^{1+2d}$ and $r\in{(0,1]}$,
\begin{equation}
\inf_{z_0\in \mathcal  G_{z,r}(K)}
 \Pp_{z_0}\left(T_{\mathcal  G_{z,r}(B)}<\tau_{\mathcal  G_{z,r}(O)}\right)
 \geqslant \ell_0,
 \label{eq:uniform-link-hitting}
\end{equation}
and
\begin{equation}
\inf_{z_0\in \mathcal  G_{z,r}(K)}
 \Pp_{z_0}\left(
 \int_0^{\tau_{\mathcal  G_{z,r}(O)}}
 \1_{\mathcal  G_{z,r}(B)}(\mathbf Z_s^{\rm aug})\dd s\geqslant\kappa_*r^2
 \right)\geqslant \ell_0.
 \label{eq:stay}
\end{equation}
\end{lemma}

\begin{proof}
Fix $z\in\R^{1+2d}$ and $r\in{(0,1]}$.  Let
$z_0=(t_0,x_0,v_0)\in \mathcal  G_{z,r}(K)$ and set 
$$
\bar z_0=(\bar t_0,\bar x_0,\bar v_0):=\mathcal  G_{z,r}^{-1}(z_0)\in K.
$$ 
{Choose the $h$ and the centered velocity control $\phi$ provided by
Definition \ref{def:directed-link-triple} for $\bar z_0$.}  Define
\[
 \phi_{r}(s):=r\phi(s/r^2),\ \ s\in[0,r^2(h+\kappa_*)].
\]
Recall that
$$
 \boldsymbol\gamma^{\bar z_0,\phi}_s
 :=\left(\bar t_0+s,
 \bar x_0+s\bar v_0+\int_0^s\phi(u)\dd u,\bar v_0+\phi(s)\right),
 $$
{A direct calculation using the definition of $\mathcal G_{z,r}$ gives}
$$
\mathcal  G_{z,r}(\boldsymbol\gamma^{\bar z_0,\phi}_{s/r^2})=\boldsymbol\gamma^{z_0,\phi_r}_{s}=
\left(t_0+s,
 x_0+sv_0+\int_0^s\phi_r(u)\dd u,v_0+\phi_r(s)\right).
$$
Moreover, {by a change of variables,} we have
\[
 \int_0^{r^2(h+\kappa_*)}|\dot\phi_{r}(s)|^2\dd s
 =
 \int_0^{h+\kappa_*}|\dot\phi(s)|^2\dd s
 \leqslant A_*.
\]
Since $d_{\rm K}(\mathcal  G_{z,r}z_1,\mathcal  G_{z,r}z_2)=r\,d_{\rm K}(z_1,z_2)$, we also have
\begin{align}\label{PP2}
d_{\rm K}\!\left(
\boldsymbol\gamma^{z_0,\phi_r}_{[0,r^2(h+\kappa_*)]},(\mathcal  G_{z,r}(O))^c
 \right)
 =rd_{\rm K}\!\left(
 \boldsymbol\gamma^{\bar z_0,\phi}_{[0,h+\kappa_*]},O^c
 \right)
 {\geqslant}2r\eta_*,
 \end{align}
 and
 \begin{align}\label{PP0}
d_{\rm K}\!\left(
 \boldsymbol\gamma^{z_0,\phi_r}_{[r^2h,r^2(h+\kappa_*)]},(\mathcal  G_{z,r}(B))^c
 \right)
 =rd_{\rm K}\!\left(
\boldsymbol\gamma^{\bar z_0,\phi}_{[h,h+\kappa_*]},B^c
 \right){\geqslant}2r\eta_*.
\end{align}
{Choose $\rho>0$, depending only on $\eta_*,h_*$, and $\kappa_*$, so that}
\begin{align}\label{PP1}
 ((h_*+\kappa_*)\rho)^{1/3}+\rho<\eta_*.
\end{align}
Define
$$
\Gamma_r:=\left\{\sup_{0\leqslant s\leqslant r^2(h+\kappa_*)}
 |\mathsf V_s-v_0-\phi_r(s)|<r\rho\right\}
 $$
 and
$$
{\widehat{\boldsymbol\beta}_u:=r\boldsymbol\beta_{r^2u},\qquad
\widehat{\mathbf A}_u:=\mathbf A_{r^2u},\qquad
\widehat W_u:=r^{-1}W_{r^2u}.}
$$
{Then $\widehat W$ is a Brownian motion for the rescaled filtration,
$|\widehat{\boldsymbol\beta}_u|\leqslant r\hbar\leqslant\hbar$, and
$\lambda I_d\preceq\widehat{\mathbf A}_u\preceq\Lambda I_d$.}
By Lemma \ref{lem:tube}, there is a constant $\ell_0=\ell_0(\Theta,\eta_*,A_*,h_*,\kappa_*)>0$
such that
\begin{align}\label{PP5}
\Pp_{z_0}\left(\Gamma_r
 \right)=\Pp\bigg(
 \sup_{s\leqslant h+\kappa_*}\left|
 \int_0^s\widehat{\boldsymbol\beta}_u\dd u
 +r\int_0^s\sqrt{2\widehat{\mathbf A}_{r^2u}}\dd \widehat W_u-\phi(s)\right|<\rho
 \bigg)
 \geqslant \ell_0.
\end{align}
{On $\Gamma_r$, we have, for $0\leqslant s\leqslant r^2(h+\kappa_*)$,}
\[
 \left|\mathsf X_s-x_0-sv_0-\int_0^s\phi_r(u)\dd u\right|
 {\leqslant}\int_0^s|\mathsf V_u-v_0-\phi_r(u)|\dd u
 \leqslant r^3(h+\kappa_*)\rho\leqslant r^3(h_*+\kappa_*)\rho,
\]
by \eqref{PP1} and the definition of $d_{\rm K}$, we have
\[
 d_{\rm K}\left(\mathbf Z_s^{\rm aug},\boldsymbol\gamma^{z_0,\phi_r}_s\right)
 <r((h_*+\kappa_*)\rho)^{1/3}\leqslant r\eta_*,
 \quad 0\leqslant s\leqslant r^2(h+\kappa_*),
\]
which together with \eqref{PP2} and \eqref{PP0} yields
$$
d_{\rm K}\left(\mathbf Z_{[0,r^2(h+\kappa_*)]}^{\rm aug},(\mathcal  G_{z,r}(O))^c\right)\geqslant r\eta_*,\ \ d_{\rm K}\left(\mathbf Z_{[r^2h,r^2(h+\kappa_*)]}^{\rm aug},(\mathcal  G_{z,r}(B))^c\right)\geqslant r\eta_*.
$$
This means that 
$$
\mathbf Z_{[0,r^2(h+\kappa_*)]}^{\rm aug}\subset  \mathcal  G_{z,r}(O),\ \ \mathbf Z_{[r^2h,r^2(h+\kappa_*)]}^{\rm aug}\subset \mathcal  G_{z,r}(B).
$$
Consequently, on the event $\Gamma_r$,
$$
T_{\mathcal  G_{z,r}(B)}<\tau_{\mathcal  G_{z,r}(O)},
$$ 
and
$$
\int_0^{\tau_{\mathcal  G_{z,r}(O)}}
 \1_{\mathcal  G_{z,r}(B)}(\mathbf Z_s^{\rm aug})\dd s
 \geqslant\int_{r^2h}^{r^2(h+\kappa_*)}
 \1_{\mathcal  G_{z,r}(B)}(\mathbf Z_s^{\rm aug})\dd s
 =\kappa_*r^2.
 $$
By \eqref{PP5}, {this completes the proof.}
\end{proof}

\begin{lemma}[{Hitting estimates for a measurable set}]
\label{lem:nfh}
Let $(K,B,O)\in\mathcal A_{\rm link}(\eta_*,A_*, h_*,\kappa_*)$ with $\kappa_*>0$.
Let $z\in \R^{1+2d}$ and $r\in{(0,1]}$.
Under \eqref{eq:ellipticity}, there is a constant
$\varepsilon_*\in(0,1/2)$ only depending on $\Theta,\eta_*,A_*,h_*,\kappa_*,H_0,|B|$, where
$H_0$ is from \eqref{eq:fixed-shape-time-span}, 
such that for any
Borel set $E\subset \mathcal  G_{z,r}(B)$ with $ |\mathcal  G_{z,r}(B)\setminus E|\leqslant\varepsilon_*|\mathcal  G_{z,r}(B)|,$
\begin{equation}
 \inf_{z_0\in \mathcal  G_{z,r}(K)}
 \Pp_{z_0}\left(T_E<\tau_{\mathcal  G_{z,r}(O)}\right)\geqslant \ell_0/2,
 \label{eq:nfh}
\end{equation}
where $\ell_0$ is the same as in Lemma \ref{lem:uniform-link-hitting}.
\end{lemma}

\begin{proof}
Fix $z_0\in \mathcal  G_{z,r}(K)$ and let $\mathbf Z_s^{\rm aug}$ be the
augmented process starting from $z_0$. {Define}
\[
 \Gamma_{z,r}:=\left\{
 \int_0^{\tau_{\mathcal  G_{z,r}(O)}}
 \1_{\mathcal  G_{z,r}(B)}(\mathbf Z_s^{\rm aug})\dd s\geqslant\kappa_*r^2
 \right\}\cap\{T_E\geqslant\tau_{\mathcal  G_{z,r}(O)}\}.
\]
Clearly,
\[
 \kappa_*r^2\1_{\Gamma_{z,r}}
 {\leqslant}
 \int_0^{\tau_{\mathcal  G_{z,r}(O)}}
 \1_{\mathcal  G_{z,r}(B)\setminus E}(\mathbf Z_s^{\rm aug})\dd s.
\]
Taking expectations and applying Theorem
\ref{thm:stopped-occupation}, we obtain
\begin{align*}
 \kappa_* r^2\Pp_{z_0}(\Gamma_{z,r})
&\leqslant\E_{z_0}\left(\int_0^{\tau_{\mathcal  G_{z,r}(O)}}
    \1_{\mathcal  G_{z,r}(B)\setminus E}(\mathbf Z_t^{\rm aug})\dd t\right)
    \leqslant C_0 r^2\left(
 \frac{|\mathcal  G_{z,r}(B)\setminus E|}{r^{4d+2}}
 \right)^{\frac1p}.
\end{align*}
Since $ |\mathcal  G_{z,r}(B)\setminus E|\leqslant\varepsilon_*|\mathcal  G_{z,r}(B)|=\varepsilon_* r^{4d+2}|B|,$ it follows from
\eqref{eq:stay} that
\begin{align*}
 \Pp_{z_0}(T_E<\tau_{\mathcal  G_{z,r}(O)})
 &\geqslant \Pp_{z_0}\left\{
 \int_0^{\tau_{\mathcal  G_{z,r}(O)}}
 \1_{\mathcal  G_{z,r}(B)}(\mathbf Z_s^{\rm aug})\dd s\geqslant\kappa_*r^2
 \right\}
   -\Pp_{z_0}(\Gamma_{z,r})\\
& \geqslant \ell_0-C_0\left(\varepsilon_*|B|\right)^{\frac1p}/\kappa_*.
\end{align*}
In particular, we may choose
$\varepsilon_*\in(0,1/2)$ small enough so that
\[
C_0\left(\varepsilon_*|B|\right)^{\frac1p} \leqslant \ell_0/2.
\]
This completes the proof.
\end{proof}

\section{Existence of weak solutions for kinetic SDEs with rough coefficients}
\label{sec:SDE}
In this section we consider the kinetic SDE
\begin{equation}
\dd\mathsf X_t=\mathsf V_t\dd t,
 \quad \dd\mathsf V_t=b(t,\mathsf X_t,\mathsf V_t)\dd t
 +\sqrt{2a(t,\mathsf X_t,\mathsf V_t)}\dd W_t,
 \quad (\mathsf X_0,\mathsf V_0)=(x_0,v_0),
 \label{eq:sde}
\end{equation}
{where
$(a,b):\R_+\times\R^d\times\R^d\to\Ssym^d\times\R^d$ are Borel
measurable and satisfy}
\begin{equation}\label{BB1}
 |b(t,x,v)|\leqslant \hbar,\quad
 \lambda I_d\preceq a(t,x,v)\preceq\Lambda I_d,
 \quad 0<\lambda\leqslant\Lambda<\infty.
\end{equation}
{We prove the existence of weak solutions to the SDE \eqref{eq:sde} with measurable diffusion coefficients.}
\begin{theorem}\label{SDE}
Under \eqref{BB1}, for every $(x_0,v_0)\in\R^{2d}$, the SDE
\eqref{eq:sde} has a weak solution starting from $(x_0,v_0)$. More precisely, there {exist} a filtered probability space 
$(\Omega,\mathcal F,(\mathcal F_t)_{t{\geqslant}0},\Pp)$ and 
$\mathcal F_t$-adapted processes $(\mathsf X_t,\mathsf V_t,W_t)$ {such} that
$W$ is an $\mathcal F_t$-Brownian motion, and{,} for each $t>0$, {$\Pp$-almost surely},
$$
\left\{
\begin{aligned}
\mathsf X_t&=x_0+tv_0+\int^t_0\mathsf V_s\dd s,\\
\mathsf V_t&=v_0+\int^t_0b(s,\mathsf X_s,\mathsf V_s)\dd s
 +\int^t_0\sqrt{2a(s,\mathsf X_s,\mathsf V_s)}\dd W_s.
 \end{aligned}
 \right.
$$
\end{theorem}

\begin{proof}
Extend $a,b$ to $\R\times\R^{2d}$ by setting
$a(t,x,v):=\lambda I_d$ for $t<0$.  Let $\rho_n$ be standard nonnegative
space--time mollifiers and set $a_n:=a*\rho_n$ and
$b_n:=b*\rho_n$.  Then
\begin{equation}\label{eq:diffusion-approximation}
 \lambda I_d\preceq a_n\preceq\Lambda I_d,
 \quad |b_n|\leqslant\hbar,
\end{equation}
and
\begin{equation}\label{eq:diffusion-convergence}
a_n\longrightarrow a,\quad b_n\longrightarrow b\quad\text{in }L^1_{\mathrm{loc}}
 ([0,\infty)\times\R^{2d}).
\end{equation}
Since {$a_n$ and $b_n$} are smooth with bounded derivatives, the SDE
\begin{equation}\label{eq:driftless-approximating-sde}
 \begin{cases}
  \dd\mathsf X_t^n=\mathsf V_t^n\dd t,\\
  \dd\mathsf V_t^n
   =b_n(t,\mathsf X_t^n,\mathsf V_t^n)+\sqrt{2a_n(t,\mathsf X_t^n,\mathsf V_t^n)}\dd W_t,
 \end{cases}
 \quad
 (\mathsf X_0^n,\mathsf V_0^n)=(x_0,v_0),
\end{equation}
has a unique nonexplosive strong solution.

Fix $T>0$. By  the Burkholder--Davis--Gundy inequality, 
for $0\leqslant s\leqslant t\leqslant T$, we have
\begin{equation}\label{eq:driftless-moment-bounds}
 \sup_n\E\sup_{t\leqslant T}|\mathsf V_t^n|^4<\infty,
 \quad
 \sup_n\E|\mathsf V_t^n-\mathsf V_s^n|^4
 \leqslant C_T|t-s|^2,
\end{equation}
and
\[
 \E|\mathsf X_t^n-\mathsf X_s^n|^4
 \leqslant |t-s|^3\int_s^t\E|\mathsf V_r^n|^4\dd r
 \leqslant C_T|t-s|^4.
\]
Hence, the laws ${\mathbf P}_n:=\Law(\mathsf X^n,\mathsf V^n, W)$ are tight on
$C([0,T];\R^{2d}\times\mR^d)$.  Since these estimates hold for every $T>0$, {a diagonal
argument yields a subsequence along which we may assume that}
\begin{equation}\label{eq:driftless-path-law-convergence}
 \mathbf P_n\Longrightarrow \mathbf P
 \quad\text{on }C([0,\infty);\R^{2d}\times\R^d)
\end{equation}
with the topology of locally uniform convergence. 
By Skorohod's representation theorem, there is a probability space 
$(\widetilde\Omega,\widetilde{\mathcal F},\widetilde{\mathbb P})$ and
processes $(\widetilde{\mathsf X}^n,\widetilde{\mathsf V}^n, \widetilde W^n)$
and $(\widetilde{\mathsf X},\widetilde{\mathsf V}, \widetilde W)$
defined on this space so that
\begin{align}\label{Sa1}
(\widetilde{\mathsf Z}^n, \widetilde W^n):=
(\widetilde{\mathsf X}^n,\widetilde{\mathsf V}^n, \widetilde W^n)\to 
(\widetilde{\mathsf X},\widetilde{\mathsf V}, \widetilde W)
=:(\widetilde{\mathsf Z}, \widetilde W),\quad\widetilde{\mathbb P}-a.s.,
\end{align}
and
$$
\mathbf P_n=\widetilde{\mathbb P}\circ(\widetilde{\mathsf X}^n,\widetilde{\mathsf V}^n, \widetilde W^n)^{-1},\quad
\mathbf P=\widetilde{\mathbb P}\circ(\widetilde{\mathsf X},\widetilde{\mathsf V}, \widetilde W)^{-1}
$$
Let
$$
\widetilde {\mathcal F}_t:=\sigma(\widetilde{\mathsf X}^n_s,\widetilde{\mathsf V}^n_s, \widetilde W^n_s, 
\widetilde{\mathsf X}_s,\widetilde{\mathsf V}_s, \widetilde W_s, s\leqslant t, n\in\mathbb N).
$$
Then $\widetilde W^n, \widetilde W$ are two $\widetilde {\mathcal F}_t$-Brownian motions, and
\begin{align}\label{Sa6}
\left\{
\begin{aligned}
\widetilde{\mathsf X}^n_t&=x_0+tv_0+\int^t_0\widetilde{\mathsf V}^n_s\dd s,\\
\widetilde{\mathsf V}^n_t&=v_0+\int^t_0b_n(s,\widetilde{\mathsf Z}^n_s)\dd s
 +\int^t_0\sqrt{2a_n(s,\widetilde{\mathsf Z}^n_s)}\dd{\widetilde W}^n_s.
 \end{aligned}
 \right.
\end{align}
Let $\widetilde{\mathbb E}$ be the expectation with respect to $\widetilde{\mathbb P}$.
Then, by Theorem~\ref{thm:occupation}, for every $f\in L^p([0,T]\times{\R^{2d}})$,
\begin{equation}\label{eq:driftless-uniform-occupation}
 \sup_n\widetilde\E\left(\int_0^T
f(t,\widetilde{\mathsf X}_t^n,\widetilde{\mathsf V}_t^n)\dd t\right)
 \leqslant  C_0\|f\|_{L^p([0,T]\times{\R^{2d}})},
\end{equation}
and also
\begin{equation}\label{eq:driftless-limit-occupation}
\widetilde\E\int_0^T
  f(t,\widetilde{\mathsf X}_t,\widetilde{\mathsf V}_t)\dd t
 \leqslant C_0\|f\|_{L^p([0,T]\times\R^{2d})}.
\end{equation}
It is well known that
\begin{align}\label{Sa7}
\lim_{n\to\infty}\widetilde{\mathbb E}\left|\int^t_0\sqrt{2a(s,\widetilde{\mathsf Z}_s)}\dd{\widetilde W}^n_s-
\int^t_0\sqrt{2a(s,\widetilde{\mathsf Z}_s)}\dd{\widetilde W}_s\right|^2=0.
\end{align}
Let $R\geqslant 1$. By \eqref{eq:driftless-uniform-occupation} with $f=|\sqrt{a_m}-\sqrt{a}|^2{\bf 1}_{B_R}$, we have
\begin{align}\label{Sa2}
&\widetilde\mE\left|\int^t_0{\bf 1}_{\{|\widetilde{\mathsf Z}^n_s|\leqslant R\}}\sqrt{2a_m(s,\widetilde{\mathsf Z}^n_s)}\dd{\widetilde W}^n_s-
\int^t_0{\bf 1}_{\{|\widetilde{\mathsf Z}^n_s|\leqslant R\}}\sqrt{2a(s,\widetilde{\mathsf Z}^n_s)}\dd{\widetilde W}^n_s\right|^2\nonumber\\
&\quad=2\widetilde\mE\int^t_0{\bf 1}_{\{|\widetilde{\mathsf Z}^n_s|\leqslant R\}}\left|\sqrt{a_m(s,\widetilde{\mathsf X}^n_s,\widetilde{\mathsf V}^n_s)}-
\sqrt{2a(s,\widetilde{\mathsf Z}^n_s)}\right|^2_{\rm F}\dd s\nonumber\\
&\quad\leqslant C_0\|\sqrt{a_m}-\sqrt{a}\|^2_{L^{2p}([0,t]\times B_R)}\longrightarrow 0,\quad\mbox{ as $m\to\infty$}.
\end{align}
On the other hand, we also have
\begin{align}\label{Sa3}
&\widetilde\mE\left|\int^t_0{\bf 1}_{\{|\widetilde{\mathsf Z}^n_s|\leqslant R\}}\sqrt{2a_m(s,\widetilde{\mathsf Z}^n_s)}\dd{\widetilde W}^n_s\right|^2\leqslant 2\Lambda d\int^t_0\widetilde{\mathbb P}\{|\widetilde{\mathsf Z}^n_s|\leqslant R\}\dd s\nonumber\\
&\qquad\qquad\leqslant \frac{2\Lambda d}{R^2}\int^t_0\widetilde{\mathbb E}|\widetilde{\mathsf Z}^n_s|^2\dd s
\stackrel{\eqref{eq:driftless-moment-bounds}}{\leqslant} C/R^2,
\end{align}
where $C$ is independent of $n$. Combining \eqref{Sa2} and \eqref{Sa3}, we get
$$
\lim_{m\to\infty}\sup_n\widetilde\mE\left|\int^t_0\sqrt{2a_m(s,\widetilde{\mathsf Z}^n_s)}\dd{\widetilde W}^n_s-
\int^t_0\sqrt{2a(s,\widetilde{\mathsf Z}^n_s)}\dd{\widetilde W}^n_s\right|^2=0,
$$
which{,} together with \eqref{Sa7}{,} yields
\begin{align}\label{Sa4}
\lim_{n\to\infty}\widetilde\mE\left|\int^t_0\sqrt{2a_n(s,\widetilde{\mathsf Z}^n_s)}\dd{\widetilde W}^n_s-
\int^t_0\sqrt{2a(s,\widetilde{\mathsf Z}_s)}\dd{\widetilde W}_s\right|^2=0.
\end{align}
Similarly, one can show
\begin{align}\label{Sa5}
\lim_{n\to\infty}\widetilde\mE\left|\int^t_0b_n(s,\widetilde{\mathsf Z}^n_s)\dd s-
\int^t_0b(s,\widetilde{\mathsf Z}_s)\dd s\right|^2=0.
\end{align}
By \eqref{Sa1}, \eqref{Sa4}, {and} \eqref{Sa5}, {passing to the limit in} \eqref{Sa6} {gives}
$$
\left\{
\begin{aligned}
\widetilde{\mathsf X}_t&=x_0+tv_0+\int^t_0\widetilde{\mathsf V}_s\dd s,\\
\widetilde{\mathsf V}_t&=v_0+\int^t_0b(s,\widetilde{\mathsf Z}_s)\dd s
 +\int^t_0\sqrt{2a(s,\widetilde{\mathsf Z}_s)}\dd{\widetilde W}_s.
 \end{aligned}
 \right.
$$
The proof is complete.
\end{proof}

{{We now fix coefficients $a,b$ satisfying
\eqref{AB2}.}  Given $z_0=(t_0,x_0,v_0)\in\R^{1+2d}$, define the shifted
coefficients}
\[
 {a^{t_0}(s,x,v):=a(t_0+s,x,v),\qquad
 b^{t_0}(s,x,v):=b(t_0+s,x,v),\qquad s\geqslant0.}
\]
{Let $(\mathsf X,\mathsf V,W)$ be any weak solution associated with
these shifted coefficients and initial condition
$(\mathsf X_0,\mathsf V_0)=(x_0,v_0)$, whose existence follows from
Theorem~\ref{SDE}, and define}
\[
 {\mathbf Z_s^{\rm aug}:=(t_0+s,\mathsf X_s,\mathsf V_s),
 \qquad \mathbf Z_0^{\rm aug}=z_0.}
\]
{We write $\Pp_{z_0}$ and $\E_{z_0}$ for the probability and
expectation corresponding to this fixed weak solution.  For a Borel set $E$
and an open set $O$, define}
\[
 {T_E:=\inf\{s\geqslant0:\mathbf Z_s^{\rm aug}\in E\},
 \qquad
 \tau_O:=\inf\{s\geqslant0:\mathbf Z_s^{\rm aug}\notin O\},}
\]
{with the convention $\inf\varnothing=+\infty$.}

\begin{lemma}
\label{lem:localized-sampling}
{Let $O\subset\R^{1+2d}$ be open, let $E\Subset O$ be a nonempty
compact set, and fix $z_0\in O$.  For every augmented weak solution
$\mathbf Z^{\rm aug}$ associated with $\cL_{a,b}$ and starting from $z_0$,
and for every nonnegative supersolution $w$ of $\cL_{a,b}$ in $O$,}
\begin{equation}
 {w(z_0)\geqslant\Pp_{z_0}(T_E<\tau_O)\inf_Ew.}
 \label{eq:localized-hitting}
\end{equation}
{In particular, \eqref{eq:localized-hitting} does not require $w$ to
extend continuously to $\partial O$.}
\end{lemma}

\begin{proof}
{Set}
\[
 {O_n:=\{z\in O:|z|<n,\ \dist(z,O^c)>n^{-1}\},
 \qquad \tau_n:=\tau_{O_n}.}
\]
{Since $E\Subset O$ and $z_0\in O$, for all sufficiently large $n$,}
\[
 {z_0\in O_n,\qquad E\Subset O_n.}
\]
{Put $\sigma_n:=T_E\wedge\tau_n$.  The generalized It\^o formula
applies to $w$ up to $\sigma_n$.  Indeed, one may mollify $w$ on a
neighborhood of $\overline O_n$ and pass to the limit using the Krylov estimate, Theorem~\ref{thm:occupation}.  Moreover, since the time component of
$\mathbf Z^{\rm aug}$ increases at unit speed and $O_n\subset\{|t|<n\}$,
$\tau_n\leqslant n-t_0<2n$ for all sufficiently large $n$.  Hence
$\sigma_n$ is bounded, and the generalized It\^o formula gives}
\begin{align*}
{\E_{z_0}w\left(\mathbf Z_{\sigma_n}^{\rm aug}\right)
 =w(z_0)+\E_{z_0}\int_0^{\sigma_n}
 \cL_{a,b}w(\mathbf Z_s^{\rm aug})\dd s.}
\end{align*}
{Since $\cL_{a,b}w\leqslant0$ Lebesgue-a.e. in $O$, $w\geqslant0$,
and $\mathbf Z_{T_E}^{\rm aug}\in E$ on $\{T_E<\tau_n\}$, we obtain}
\begin{align*}
 {w(z_0)}
 {\geqslant\E_{z_0}w(\mathbf Z_{\sigma_n}^{\rm aug})}
 {\geqslant\E_{z_0}\!\left[
 w(\mathbf Z_{T_E}^{\rm aug})\1_{\{T_E<\tau_n\}}\right]}
 {\geqslant\Pp_{z_0}(T_E<\tau_n)\inf_Ew.}
\end{align*}
{On the event $\{T_E<\tau_O\}$, continuity implies that the trajectory
up to time $T_E$ is compact and lies at a positive distance from $O^c$.
Since $O_n\uparrow O$, it follows that}
\[
 {\{T_E<\tau_n\}\uparrow\{T_E<\tau_O\}.}
\]
Letting $n\to\infty$ in the preceding inequality, we get
\eqref{eq:localized-hitting}.  This completes the proof.
\end{proof}

\section{Harnack inequalities of kinetic equations in nondivergence form}

{Throughout this section, we consider the following second-order kinetic
operator:}
\[
 \cL_{a,b}:=\partial_t+v\cdot\nabla_x
+a(t,x,v):\nabla_v^2 +b(t,x,v)\cdot\nabla_v,
\]
{where $a:\R^{1+2d}\to\Ssym^d$ and
$b:\R^{1+2d}\to\R^d$ are Borel measurable and satisfy}
\begin{align}\label{AB2}
 \lambda I_d\preceq a\preceq\Lambda I_d,
 \qquad |b|\leqslant \hbar.
\end{align}
\begin{definition}\label{Def0}
{Let $\Omega\subset\R^{1+2d}$ be open.  A function
$w\in C(\Omega)$ whose distributional derivatives satisfy}
\[
 {\partial_t w,\ \nabla_xw,\ \nabla_vw,\ \nabla_v^2w
 \in \cap_{p>1}L^p_{\mathrm{loc}}(\Omega)}
\]
is called a solution (respectively, a supersolution) of
$\cL_{a,b}$ in $\Omega$ if
\[
 {\cL_{a,b}w=0
 \quad\text{(respectively, $\cL_{a,b}w\leqslant0$)}
 \quad\text{Lebesgue-a.e. in }\Omega.}
\]
{If, in addition, $w\geqslant0$ in $\Omega$, it is called a
nonnegative solution (respectively, a nonnegative supersolution).}
\end{definition}
{The following invariance property is a direct consequence of the chain
rule.}
\begin{lemma}[{Invariance under kinetic translations and dilations}]\label{lem_scaling_PDE}
{Let $w$ be a solution or a supersolution of $\cL_{a,b}w=0$ in
$\Omega$.}
{Then, for every
$(\widetilde z,\widetilde r)\in\R^{1+2d}\times(0,\infty)$,}
\[
 {\widetilde w(z):=w(\mathcal G_{\widetilde z,\widetilde r}(z))}
\]
{is, respectively, a solution or a supersolution of
$\cL_{\widetilde a,\widetilde b}$ in
$\mathcal G_{\widetilde z,\widetilde r}^{-1}(\Omega)$, where}
\[
 {\widetilde a(z)=a(\mathcal G_{\widetilde z,\widetilde r}(z)),
 \qquad
 \widetilde b(z)=\widetilde r\,
 b(\mathcal G_{\widetilde z,\widetilde r}(z)).}
\]
{More precisely,}
\[
 {\cL_{\widetilde a,\widetilde b}\widetilde w(z)
 =\widetilde r^2
 (\cL_{a,b}w)(\mathcal G_{\widetilde z,\widetilde r}(z)).}
\]
\end{lemma}

{We introduce three reference sets.  Let $Q_+$ and $Q_-$ be the closed
cylinders}
\begin{align*}
Q_+&:=\{(t,x,v):-1\leqslant t\leqslant 0,\ |x|\leqslant 1,\ |v|\leqslant 1\},\\
Q_-&:=\{(t,x,v): -10\leqslant t\leqslant-9,\ |x|\leqslant 1,\ |v|\leqslant 1\},
\end{align*}
{and let $Q_0$ be the open product cylinder}
\[
 Q_0:=\{(t,x,v): |t|<1,\ |x|<1,\ |v|<1\}.
\]
{In what follows, we use the parameter tuple}
\[
 \Theta:=(d,\lambda,\Lambda,\hbar).
\]
\subsection{Weak Harnack inequality}
{In this subsection, we establish the weak Harnack inequality for
nonnegative supersolutions of $\cL_{a,b}$.  The following double-ball and
critical-density estimates follow directly from Lemmas
\ref{lem:localized-sampling}, \ref{lem:uniform-link-hitting}, and
\ref{lem:nfh}.}
\begin{lemma}[{Double-ball and critical-density estimates}]
\label{lem:double-ball}
Let $(K,B,O)\in\mathcal A_{\rm link}(\eta_*,A_*, h_*,\kappa_*)$.
{Under \eqref{AB2}, there exists
$\ell_0=\ell_0(\Theta,\eta_*,A_*,h_*,\kappa_*)\in(0,1)$ such that, for every
$z\in\R^{1+2d}$, every $r\in(0,1]$, and every nonnegative supersolution
$w$ of $\cL_{a,b}$ in $\mathcal G_{z,r}(O)$,}
\begin{equation}
 \inf_{\mathcal G_{z,r}(K)}w\geqslant \ell_0\inf_{\mathcal G_{z,r}(B)}w.
 \label{eq:double-ball}
\end{equation}
{If $\kappa_*>0$, there also exist $\mu,\ell_1\in(0,1/2)$, depending
only on $\Theta,\eta_*,A_*,h_*,\kappa_*,H_0$, and $|B|$, where $H_0$ is
defined in \eqref{eq:fixed-shape-time-span}, such that, for the same
$z,r,w$ and every $L>0$, if}
\begin{equation}
 |\mathcal G_{z,r}(B)\setminus\{w\geqslant L\}|
 \leqslant \mu|\mathcal G_{z,r}(B)|,
 \label{eq:cd-assumption}
\end{equation}
{then}
\begin{equation}
 \inf_{\mathcal G_{z,r}(K)}w\geqslant \ell_1 L,
 \quad{\text{or, equivalently,}}\quad
 \mathcal G_{z,r}(K)\subset\{w\geqslant\ell_1 L\}.
 \label{eq:critical-density}
\end{equation}
\end{lemma}

\begin{proof}
Fix $z_0\in\mathcal G_{z,r}(K)$.
{Let $\mathbf Z_s^{\rm aug}$ be any augmented weak solution associated
with $\cL_{a,b}$ and starting from $z_0$.}
By Lemma \ref{lem:localized-sampling} and 
\eqref{eq:uniform-link-hitting}, we have
\[
 w(z_0)\geqslant \Pp_{z_0}\left(T_{\mathcal G_{z,r}(B)}<\tau_{\mathcal G_{z,r}(O)}\right)\inf_{\mathcal G_{z,r}(B)}w
 \geqslant \ell_0\inf_{\mathcal G_{z,r}(B)}w.
\]
{For a fixed $L>0$, set
$E:=\{w\geqslant L\}\cap\mathcal G_{z,r}(B)$. Hence by Lemmas
\ref{lem:localized-sampling} and \ref{lem:nfh}, }
\[
 w(z_0)\geqslant \Pp_{z_0}\left(T_{E}<\tau_{\mathcal G_{z,r}(O)}\right) L
 \geqslant \ell_0 L/2.
\]
{Taking $\mu:=\varepsilon_*$ from Lemma \ref{lem:nfh} and
$\ell_1:=\ell_0/2$ completes the proof.}
\end{proof}

{Before proving the weak Harnack inequality, we apply the ink-spots
theorem, Theorem~\ref{thm:ink-spots}, to obtain the following key level-set
estimate for nonnegative supersolutions.}
\begin{lemma}\label{lem:applying_ink-spot}
{For every $m\in\mathbb N$, under \eqref{AB2}, there exist constants
$\mu=\mu(\Theta)\in(0,1/2)$,
$C_m=C_m(\Theta,m)>0$, $M_m=M_m(\Theta,m)>1$, and
$R_m=R_m(d,m)>4$ such that, for every nonnegative supersolution $w$ of
$\cL_{a,b}$ in $\delta_{R_m}(Q_0)$ satisfying
$\inf_{Q_-}w<1$, and for every $M>M_m$,}
\begin{align}\label{proof:iteration_k}
|\{z:w(z)>M^k\}\cap Q_+|
\leqslant C_m\left(\frac{m+1}{m}(1-c_0\mu)\right)^k,
\quad {k=0,1,2,\ldots,}
\end{align}
{where $c_0=c_0(d)\in(0,1)$ is the constant in Theorem
\ref{thm:ink-spots}.}
\end{lemma}
\begin{proof}
Fix $m\in\mathbb N$.
Let $w$ be a nonnegative supersolution of $\cL_{a,b}$ in $\delta_{R_m}(Q_0)$ with
$\inf_{Q_-}w<1$, where $R_m>4$ will be chosen below.  {For every such
choice, we have}
\[
 Q_+\cup Q_-\Subset\delta_{R_m}(Q_0).
\]
Define
\[
 K_0:=\{-1\leqslant t\leqslant-\tfrac12,\ |x|\leqslant1,\ |v|\leqslant1\}.
\]
{By Lemma \ref{lem:elementary-directed-steering}, there exist a bounded
open set $O_0$ and constants $\eta_0,A_0,h_0,\kappa_0>0$, depending only on
$d$, such that}
\[
 (K_0,Q_+,O_0)
 \in\mathcal A_{\rm link}(\eta_0,A_0,h_0,\kappa_0).
\]
Let $\mu\in(0,1/2)$ and $\ell_{\rm cd}\in(0,1)$ be the constants in the
critical estimate of Lemma \ref{lem:double-ball}.
{Let $M>1$, to be fixed below.  For $k=0,1,2,\ldots$, define}
\[
\mathcal E^k_{\mu}:=
\Big\{(z,r)\in Q_+\times(0,1): Q_{z,r} \subset Q_+,\
|Q_{z,r}\setminus\{w>M^k\}| \leqslant \mu|Q_{z,r}|\Big\},
\]
{where $Q_{z,r}:=\mathcal G_{z,r}(Q_+)$.}
{Suppose that $(z,r)\in\mathcal E^k_\mu$ and}
\begin{align}\label{0904:00}
    {\mathcal G_{z,r}(O_0)\subset\delta_{R_m}(Q_0).}\tag{H$_1$}
\end{align}
{By applying Lemma \ref{lem:double-ball} with
$(K,B,O)=(K_0,Q_+,O_0)$ and $L=M^k$, we have}
\begin{align}\label{proof:CD}
\inf_{{\mathcal G_{z,r}(K_0)}}w
    \geqslant \ell_{\rm cd} M^k,
\quad{\text{and hence}}\quad
\mathcal G_{z,r}(K_0)\subset\{w\geqslant \ell_{\rm cd} M^k\}.
\end{align}
{We divide the remainder of the proof into four steps.}

\medskip
\noindent\textbf{Step 1:} 
{We show that if $M$ is chosen sufficiently large, depending on $m$,
then, for every $(z,r)\in\mathcal E^k_\mu$,}
\begin{align}\label{proof:H1}
{\bar Q_{z,r}^{\,m}}\subset \{w>M^{k-1}\},
\end{align}
where
$$
{\bar Q_{z,r}^{\,m}}=\{-(m+1)r^2\leqslant t-t_z<-r^2,\, |x-x_z-(t-t_z)v_z|\leqslant (m+2)r^3,\, |v-v_z|\leqslant r \}.
$$
{Set}
\[
{\bar Q^{\,m}:=\{-(m+1)\leqslant t\leqslant -1,\ |x|\leqslant m+2,\ |v|\leqslant1\}.}
\]
{By Lemma \ref{lem:elementary-directed-steering}, there exist a bounded
open set $O_m$ and constants $\eta_m,A_m,h_m,\kappa_m>0$, depending only on
$d$ and $m$, such that}
\[
({\bar Q^{\,m}},K_0,O_m)
\in\mathcal A_{\rm link}(\eta_m,A_m,h_m,\kappa_m).
\]
{Since $\bar Q_{z,r}^{\,m}\subset
\mathcal G_{z,r}(\bar Q^{\,m})$, by Lemma \ref{lem:double-ball} with
$(K,B,O)=(\bar Q^{\,m},K_0,O_m)$, there is a constant
$\ell_m=\ell_m(\Theta,m)\in(0,1)$ such that, whenever}
\begin{align}\label{0904:01}
  {\mathcal G_{z,r}(O_m)\subset\delta_{R_m}(Q_0),}\tag{H$_2$}
\end{align}
{holds (so that $w$ is a nonnegative supersolution of $\cL_{a,b}$ in
$\mathcal G_{z,r}(O_m)$),}
\begin{align*}
    \inf_{{\bar Q_{z,r}^{\,m}}} w
    \geqslant \inf_{{\mathcal G_{z,r}(\bar Q^{\,m})}} w
    \geqslant \ell_m \inf_{{\mathcal G_{z,r}(K_0)}} w.
\end{align*}
{Together with \eqref{proof:CD}, this implies}
\begin{align*}
\inf_{{\bar Q_{z,r}^{\,m}}} w\geqslant \ell_m\ell_{\rm cd}M^k.
\end{align*}
{Therefore, \eqref{proof:H1} holds whenever
$M>(\ell_m\ell_{\rm cd})^{-1}$.}

\medskip
\noindent\textbf{Step 2:}
{For $j=0,1,2,\ldots$, set}
\[
 {z_j:=(1-2^{2j},0,0),
 \qquad K_j:=\mathcal G_{z_j,2^j}(K_0).}
\]
{By Lemma \ref{lem:elementary-directed-steering}, there exist a bounded open set $O_*$ and
constants $\eta_*,A_*,h_*>0$, depending only on $d$, such that}
\[
 (K_1,K_0,O_*)\in\mathcal A_{\rm link}(\eta_*,A_*,h_*,0).
\]
{Fix $(z,r)\in\mathcal E^k_\mu$, set
$j_r:=\lceil\log_2(1/r)\rceil$, and, for $j\geqslant0$, let
$\bar z_j:=\mathcal G_{z,r}(z_j)$.  Since
$\mathcal G_{z_j,2^j}(z_\ell)=z_{j+\ell}$ for $j,\ell\geqslant0$,
by \eqref{KiDi2}, we have}
\begin{align}\label{0906:00}
 \mathcal G_{z,r}(K_{j+1}){=\mathcal G_{\bar z_j,r\,2^j}(K_1),}\quad
 \mathcal G_{z,r}(K_j){=\mathcal G_{\bar z_j,r\,2^j}(K_0).}
\end{align}
{For $0\leqslant j\leqslant j_r-1$, one has $r\,2^j<1$.
Consequently, by Lemma \ref{lem:double-ball} with
$(K,B,O)=(K_1,K_0,O_*)$, there is a constant
$\ell_*=\ell_*(\Theta)\in(0,1)$ such that}
\[
 {\inf_{\mathcal G_{z,r}(K_{j+1})}w
 \geqslant\ell_*
 \inf_{\mathcal G_{z,r}(K_j)}w,
 \qquad 0\leqslant j\leqslant j_r-1,}
\]
{provided that}
\begin{equation}\label{0904:02}
 {\mathcal G_{\bar z_j,r\,2^j}(O_*)
 \subset\delta_{R_m}(Q_0),
 \qquad 0\leqslant j\leqslant j_r-1.}\tag{H$_3$}
\end{equation}
{Iterating these inequalities and using \eqref{proof:CD}, we obtain}
\begin{equation}\label{0831:00}
 {\inf_{\mathcal G_{z,r}(K_{j_r})}w
 \geqslant\ell_*^{j_r}\ell_{\rm cd}M^k.}
\end{equation}
 {Choose
$z^\dagger_0=(-2/3,0,0)\in K_0$
and set}
\[
 \rho_0:=d_{\rm K}(z^\dagger_0,K_0^c)=1/\sqrt6>0.
\]
Set
\[
 {z^\dagger_r
 :=\mathcal G_{\mathcal G_{z,r}(z_{j_r}),r\,2^{j_r}}(z^\dagger_0)=\mathcal G_{\bar z_{j_r},\,r\,2^{j_r}}(z^\dagger_0)
 =:(t^\dagger_r,x^\dagger_r,v^\dagger_r).}
\]
{By the definition of $j_r$, $r\,2^{j_r}\in[1,2)$.} Thus, by \eqref{0906:00} and \eqref{KiDi1}, we have
\[
 {d_{\rm K}(z^\dagger_r,\mathcal G_{z,r}(K_{j_r})^c)=d_{\rm K}(\mathcal G_{\bar z_{j_r},\,r\,2^{j_r}}(z^\dagger_0)
, \mathcal G_{\bar z_{j_r},\,r\,2^{j_r}}(K_0)^c
)
 =(r\,2^{j_r})\rho_0\geqslant\rho_0.}
\]
{The time coordinate of $z^\dagger_r$ satisfies}
\begin{align*}
 t^\dagger_r
 {=t_z+r^2(1-2^{2j_r})-\tfrac23(r\,2^{j_r})^2}{\geqslant-1-\tfrac53(r\,2^{j_r})^2>-8
 =1+\sup\{t:(t,x,v)\in Q_-\}.}
\end{align*}
{Furthermore,}  
{since by \eqref{KiDi2},} 
\[
 {\mathcal G_{z,r}(K_{j_r})
 =\left\{z\circ\delta_r(z_{j_r})
 \circ\delta_{r\,2^{j_r}}(z_0):z_0\in K_0\right\},}
\]
{it follows from \eqref{Af3} that}
\begin{align*}
 {\sup_{z'\in \mathcal G_{z,r}(K_{j_r})}\normK{z'}}
 &{\leqslant
 \sup_{z\in Q_+}\normK{z}
 +\normK{\delta_r(z_{j_r})}
 +\sup_{z_0\in K_0}\normK{\delta_{r\,2^{j_r}}(z_0)}}{\leqslant1+2(r\,2^{j_r})\leqslant5.}
\end{align*}
{In particular, there is
a constant $C_d<\infty$, independent of $(z,r)$, such that}
\[
 {t^\dagger_r-\sup\{t:(t,x,v)\in Q_-\}>1,
 \quad d_{\rm K}(z^\dagger_r,\mathcal G_{z,r}(K_{j_r})^c)\geqslant\rho_0,
 \quad
 \sup_{z\in Q_-, z'\in \mathcal G_{z,r}(K_{j_r})}
 d_{\rm K}(z,z')\leqslant C_d.}
\]
{By Lemma
\ref{lem:elementary-directed-steering}, there are constants $\eta,A,h>0$ and a bounded open set $\widetilde O$,
independent of $(z,r)$, such that}
\[
 {(Q_-,\mathcal G_{z,r}(K_{j_r}),\widetilde O)
 \in\mathcal A_{\rm link}(\eta,A,h,0).}
\]
{Now, applying Lemma \ref{lem:double-ball} and \eqref{0831:00}, we
obtain a constant $\ell_{\rm H}\in(0,1)$, independent of $(z,r)$, such that}
\[
 \inf_{Q_-}w
 \geqslant\ell_{\rm H}\inf_{B_{z,r}}w
 \geqslant\ell_{\rm H}\ell_*^{j_r}\ell_{\rm cd}M^k,
\]
provided that
\begin{equation}\label{0904:03}
 \widetilde O\subset\delta_{R_m}(Q_0).\tag{H$_4$}
\end{equation}
{Since $j_r\leqslant\log_2(1/r)+1$, set}
\[
 \beta:=\log_2(1/\ell_*)>0,
 \qquad c_*:=\ell_{\rm H}\ell_*\ell_{\rm cd}\in(0,1).
\]
{By the assumption $\inf_{Q_-}w<1$, we have}
\[
 1>\ell_{\rm H}\ell_*^{j_r}\ell_{\rm cd}M^k
 \geqslant c_*r^\beta M^k.
\]
{Therefore, with $r_k:=(c_*M^k)^{-1/\beta}$, assumptions
\eqref{0904:00}--\eqref{0904:03} imply}
\begin{equation}\label{proof:r_k}
 {r<r_k\qquad\text{for every }(z,r)\in\mathcal E^k_\mu.}
\end{equation}

\medskip
\noindent\textbf{Step 3:}
{Set}
\[
 \theta_m:=\frac{m+1}{m}(1-c_0\mu)>0.
\]
Choose
\[
 M_m>
 \max\left\{1,(\ell_m\ell_{\rm cd})^{-1},c_*^{-1},
 \left(\frac{2}{\theta_m}\right)^{\beta/2}\right\},
\]
{and fix $M>M_m$.  By Steps 1 and 2, for every
$(z,r)\in\mathcal E^k_\mu$ with $k\geqslant1$, assumptions
\eqref{0904:00}--\eqref{0904:03} imply}
\[
 r<r_k:=(c_*M^k)^{-1/\beta}<1,
 \qquad
 {\bar Q_{z,r}^{\,m}}\subset\{w>M^{k-1}\}.
\]
{We may therefore apply Theorem \ref{thm:ink-spots} with}
\[
 E:=\{w>M^k\}\cap Q_+,\quad
 F:=\{w>M^{k-1}\},\quad r_0:=r_k,
\]
{and for $a_k:=|\{w>M^k\}\cap Q_+|$,}
\begin{align*}
 a_k
 &\leqslant \theta_m\bigl(a_{k-1}+C_0m r_k^2\bigr)=\theta_m a_{k-1}
   +C_0m\theta_m (c_*M^k)^{-2/\beta}.
\end{align*}
{Set $q:=M^{-2/\beta}$ and
$D_m:=C_0m\theta_m c_*^{-2/\beta}$.  The choice of $M_m$ guarantees
$q\leqslant\theta_m/2$.  Iterating the preceding recurrence and using
$a_0\leqslant|Q_+|$, we obtain}
\begin{align*}
 a_k
 &\leqslant \theta_m^k a_0
   +D_m\sum_{j=1}^k\theta_m^{k-j}q^j
   \leqslant \theta_m^k\left(
       |Q_+|+D_m\sum_{j=1}^{\infty}
       \left(\frac q{\theta_m}\right)^j\right)
 \leqslant (|Q_+|+D_m)\theta_m^k.
\end{align*}
{Thus \eqref{proof:iteration_k} holds with
$C_m:=|Q_+|+D_m$, which depends only on $\Theta$ and $m$.}

\medskip
\noindent\textbf{Step 4:}
{It remains to choose $R_m$ for which
\eqref{0904:00}--\eqref{0904:03} hold.  Define}
\begin{align}\label{eq:uniform-outer-union}
U_m:=\left(\cup_{(z,r)\in Q_+\times(0,1)}\mathcal G_{z,r}(O_0\cup O_m)\right)\cup\left(\cup_{z\in Q_+, r\in(0,1)}\cup_{j=0}^{j_r} {\mathcal G_{\bar z_j, r\,2^j}}(O_*)\right)\cup \tilde{O}
\end{align}
{Since $Q_+$, $O_0$, and $O_m$ are bounded, the first union in
\eqref{eq:uniform-outer-union} is bounded; more precisely,}
\begin{align*}
\sup_{z\in Q_+,0<r<1}
 \sup_{z'\in\mathcal G_{z,r}(O_0\cup O_m)}\normK{z'}
&{\leqslant1+
 \sup_{z''\in O_0\cup O_m}\normK{z''}<\infty.}
\end{align*}
{For the second union, set $s_j:=r\,2^j$.  If
$0\leqslant j\leqslant j_r-1$, then $s_j\in(0,1)$ and}
\[
 {\delta_r(z_j)=(r^2-s_j^2,0,0),
 \qquad \bar z_j=z\circ(r^2-s_j^2,0,0).}
\]
{Hence}
\begin{align*}
\sup_{z'\in\mathcal G_{\bar z_j,s_j}(O_*)}\normK{z'}
&{\leqslant\normK{z}
 +\normK{(r^2-s_j^2,0,0)}
 +\sup_{z''\in O_*}\normK{\delta_{s_j}z''}}
 {\leqslant2+\sup_{z''\in O_*}\normK{z''}<\infty.}
\end{align*}
{Thus the second union is bounded uniformly in $z,r,j$.  Since
$\widetilde O$ is bounded, $U_m$ is bounded.  Choosing}
\[
 {R_m:=5+\sup_{z\in U_m}\normK{z}<\infty,}
\]
{we have
$U_m\subset\delta_{R_m}(Q_0)$ and $R_m>4$.  This completes the proof.}
\end{proof}

{We can now prove the weak Harnack inequality.}
\begin{theorem}[Weak Harnack inequality]
\label{thm:weak-harnack}
{Under \eqref{AB2}, there exist
$m=m(\Theta)\in\mathbb N$, $p_0=p_0(\Theta)>0$, and
$C_0=C_0(\Theta)>0$ such that every nonnegative
supersolution $w$ of $\cL_{a,b}$ in $\delta_{R_m}(Q_0)$, where $R_m$ is
given by Lemma \ref{lem:applying_ink-spot}, satisfies}
\begin{equation}
 \left(\fint_{Q_+}w^{p_0}\right)^{1/p_0}
 \leqslant {C_0} \inf_{Q_-}w.
 \label{eq:weak-harnack}
\end{equation}
\end{theorem}

\begin{proof}
{Let $\mu\in(0,1/2)$ be the constant in Lemma
\ref{lem:applying_ink-spot}, and choose $m\in\mathbb N$ sufficiently large
that}
\[
\theta_m:=\frac{m+1}{m}(1-c_0\mu)\in(0,1).
\]
{Fix $M:=2M_m$, where $M_m$ is the constant in Lemma
\ref{lem:applying_ink-spot}, and define}
\begin{equation}
 p_0:=\frac{\log(1/\theta_m)}{2\log M}>0.
 \label{eq:pwh-choice}
\end{equation}
{For $\eta>0$, let}
\[
 u_\eta:=\frac{w}{\inf_{Q_-} w+\eta}.
\]
{Then $u_\eta$ is a nonnegative supersolution of
$\cL_{a,b}$ in $\delta_{R_m}(Q_0)$, and}
\begin{equation}
 \inf_{Q_-}u_\eta<1.
 \label{eq:wh-normalization}
\end{equation}
{For every $z\in Q_+$,} note that
\[
 {u_\eta^{p_0}(z)
 \leqslant1+\sum_{k=0}^\infty
 M^{p_0(k+1)}\1_{\{u_\eta>M^k\}\cap Q_+}(z).}
\]
{Hence by Lemma \ref{lem:applying_ink-spot} and \eqref{eq:pwh-choice},}
\begin{align*}
 \fint_{Q_+}u_\eta^{p_0}
 \leqslant
 1+
 {\frac{C_m}{|Q_+|}}M^{p_0}
 \sum_{k=0}^\infty
 \left(M^{p_0}\theta_m\right)^k
 =
 1+
 {\frac{C_m}{|Q_+|}}M^{p_0}
 \sum_{k=0}^\infty
 \theta_m^{k/2}
 =:C_1<\infty.
\end{align*}
{Since $w=(\inf_{Q_-}w+\eta)u_\eta$, we have}
\[
 {\left(\fint_{Q_+}w^{p_0}\right)^{1/p_0}
 \leqslant
 C_1^{1/p_0}(\inf_{Q_-}w+\eta).}
\]
{Letting $\eta\downarrow0$, we prove \eqref{eq:weak-harnack} with
$C_0:=C_1^{1/p_0}$.}
\end{proof}

\subsection{Strong Harnack inequality}
\label{sec:strong-harnack}

{In this subsection, we prove the strong Harnack inequality by
establishing a local $L^{p_0}$--$L^\infty$ estimate for nonnegative solutions
of $\cL_{a,b}$.}
\begin{theorem}[Strong Harnack inequality]
\label{thm:strong-harnack}
{Assume \eqref{AB2}, and let $m\in\mathbb N$ and $R_m>4$ be the
constants appearing in Theorem \ref{thm:weak-harnack}.  If $w$ is a
nonnegative solution of $\cL_{a,b}$ in $\delta_{R_m}(Q_0)$, then, for every
$Q\Subset\mathring Q_+$, there exists a constant}
$C=C\bigl(\Theta,d_{\rm K}(Q,Q_+^c)\bigr)>0$
{such that}
  \begin{align*}
      \sup_{Q} w\leqslant C \inf_{Q_-}w. 
  \end{align*}
\end{theorem}

\begin{proof}
{By the weak Harnack inequality \eqref{eq:weak-harnack}, it is enough
to prove the following local estimate for every nonnegative solution in an
open neighborhood of $Q_+$:}
\begin{align*}
      \sup_{Q} w\leqslant C \left(\int_{Q_+}w^{p_0}\right)^{1/p_0}.
  \end{align*}
{Without loss of generality, we assume $\sup_{Q} w>0$.
Set $N:=4d+2$ and define}
$$
 {M:=\max_{z\in Q_+}\varrho(z)^{N/p_0}w(z)\in(0,\infty),}
$$
where $\varrho(z):=d_{\rm K}(z,Q_+^c)\leqslant1$.
{Since the function $z\mapsto\varrho(z)^{N/p_0}w(z)$ is continuous on the
compact set $Q_+$, and vanishes on $\partial Q_+$, there is a point
$z_0=(t_0,x_0,v_0)\in\mathring Q_+$ such that}
\begin{align}\label{0902:01}
 {\varrho(z_0)^{N/p_0}w(z_0)=M.}
\end{align}
Recall $Q_0:=\{|t|<1,\ |x|<1,\ |v|<1\}$ and set
\[
 r_0:=\tfrac{\varrho(z_0)}{2}>0.
\]
Note that
$d_{\rm K}(0,z)< 1$ for $z\in Q_0$.  By left invariance and homogeneity of
$d_{\rm K}$, it follows that
\[
 \mathcal G_{z_0,r_0}(Q_0)\subset Q_+,
 \quad
 d_{\rm K}(z_0,z')\leqslant \tfrac{\varrho(z_0)}{2}=r_0
 \quad\text{for }z'\in\overline{\mathcal G_{z_0,r_0}(Q_0)}.
\]
For every $z'\in\overline{\mathcal G_{z_0,r_0}(Q_0)}$ and every $z\in Q_+^c$, by \eqref{KiDi1} we have
\[
 2r_0=\varrho(z_0)\leqslant d_{\rm K}(z_0,z)
 \leqslant d_{\rm K}(z_0,z')+d_{\rm K}(z',z)\leqslant {r_0}+d_{\rm K}(z',z),
\]
and hence
\[
 \varrho(z')\geqslant r_0.
\]
Consequently, by \eqref{0902:01},
\begin{align}\label{0902:03}
 w(z')
 \leqslant \frac{M}{\varrho(z')^{N/p_0}}
 \leqslant \frac{M}{r_0^{N/p_0}}
 \leqslant 2 w(z_0)
 \quad\text{for }z'\in\overline{\mathcal G_{z_0,r_0}(Q_0)}.
\end{align}
{Define}
\[
 E:=\left\{z\in\mathcal G_{z_0,r_0}(Q_0):
               w(z)>\frac{w(z_0)}2\right\}.
\]
{Let
$\mathbf Z_s^{\rm aug}:=(t_0+s,\mathsf X_s,\mathsf V_s)$ be an augmented
weak solution associated with $\cL_{a,b}$ and starting from $z_0$, and set}
\[
 {\tau_v:=\inf\{s\geqslant0:|\mathsf V_s-v_0|\geqslant r_0\},}
 \quad
 {\tau_0:=\inf\left\{s\geqslant0:
 w(\mathbf Z_s^{\rm aug})\leqslant\frac{w(z_0)}2\right\}.}
\]
We divide the proof into two steps.

\medskip
\noindent\textbf{Step 1.}
{We first prove that there exists a constant}
$c_1=c_1(\Theta)\in(0,1)$ such that
\begin{align}\label{0902:10}
 {\Pp_{z_0}(\tau_v\wedge\tau_0>c_1r_0^2)\geqslant1/6.}
\end{align}
{Put $T_0:=c_1r_0^2$.  Note that for  $s<\tau_v\wedge T_0$,}
\[
 {|\mathsf X_s-x_0-sv_0|
\leqslant\int_0^s|\mathsf V_r-v_0|\dd r}
 \leqslant r_0T_0=c_1r_0^3.
\]
{Thus, if $c_1<1$, then}
\begin{align}\label{0902:08}
 {\mathbf Z_s^{\rm aug}\in\overline{\mathcal G_{z_0,r_0}(Q_0)}}
 \quad\text{for }0\leqslant s\leqslant\tau_v\wedge T_0,
\end{align}
{and by \eqref{0902:03},}
\begin{align}\label{0902:05}
 {w(\mathbf Z_s^{\rm aug})\leqslant2w(z_0),}
 \quad 0\leqslant s\leqslant\tau_v\wedge T_0.
\end{align}
Let $\tau:=\tau_0\wedge\tau_v\wedge T_0$ and
$A:=\{\tau_0\leqslant\tau_v\wedge T_0\}$.  By the generalized stopped It\^o
formula {and \eqref{0902:05}, we have}
\begin{align*}
w(z_0)=\E_{z_0}w(\mathbf Z_\tau^{\rm aug})
 &=\E_{z_0}\Big(w(\mathbf Z_{\tau_0}^{\rm aug}){\bf 1}_{A}\Big)
 +\E_{z_0}\Big(w(\mathbf Z_{\tau_v\wedge T_0}^{\rm aug}){\bf 1}_{A^c}\Big)\\
 &\leqslant\frac{w(z_0)}2\Pp_{z_0}(A)
     +2w(z_0)\Pp_{z_0}(A^c).
\end{align*}
{Therefore}
\begin{align}\label{0902:06}
 {\Pp_{z_0}\bigl(\tau_0>\tau_v\wedge T_0\bigr)
 =\Pp_{z_0}(A^c)
 \geqslant\tfrac13.}
\end{align}
{Since $r_0\leqslant1$, if $c_1\hbar\leqslant1/2$, then  the martingale}
\[
 \mathcal M_s
{:=\mathsf V_s-v_0-\int_0^s b(\mathbf Z_r^{\rm aug})\dd r
 =\int_0^s\sqrt{2a(\mathbf Z_r^{\rm aug})}\dd W_r}
\]
satisfies{,} on $\{\tau_v\leqslant T_0\}$,
\[
\sup_{0\leqslant s\leqslant T_0}|\mathcal M_s|
 \geqslant\frac{r_0}{2}.
\]
Hence, by Chebyshev's inequality, Doob's maximal inequality, and the bound
$a\preceq\Lambda I_d$, we have
\begin{align}\label{0902:07}
 \Pp_{z_0}(\tau_v\leqslant T_0)
{\leqslant\frac{4}{r_0^2}\E_{z_0}\left(
 \sup_{0\leqslant s\leqslant T_0}|\mathcal M_s|^2\right)}{\leqslant\frac{16}{r_0^2}
 \E_{z_0}|\mathcal M_{T_0}|^2
 \leqslant32d\Lambda c_1.}
\end{align}
{Choose $c_1=c_1(\Theta)\in(0,1)$ so that
$c_1\hbar\leqslant1/2$ and $32d\Lambda c_1\leqslant 1/6$.  
Subtracting the exceptional event $\{\tau_v\leqslant T_0\}$ from
the event in \eqref{0902:06}, we get \eqref{0902:10}.}

\par\medskip
\noindent\textbf{Step 2.}
{On $\{\tau_v\wedge\tau_0>T_0\}$, we have
$\mathbf Z_s^{\rm aug}\in E$ for every $0\leqslant s\leqslant T_0$.
Hence}
\[
 {\E_{z_0}\left(\int_0^{\tau_{\mathcal G_{z_0,r_0}(Q_0)}}
       \1_E(\mathbf Z_s^{\rm aug})\dd s\right)
 \geqslant\frac{c_1r_0^2}6.}
\]
Now, we apply the stopped occupation estimate \eqref{eq:stopped-occupation} with
$K=\{(0,0,0)\}$, $O=Q_0$, $z=z_0${, and} $r=r_0$, to obtain
\[
 \frac{c_1r_0^2}6
 \leqslant  C(\Theta)r_0^2
\left(\frac{|E|}{r_0^N}\right)^{\frac12-\frac1{2\digamma}}.
\]
Since $1-1/\digamma>0$, there is a constant
$c_3=c_3(\Theta)>0$ such that
\begin{align}\label{0902:04}
 {|E|\geqslant c_3r_0^N.}
\end{align}
Thus, by the definition of $E$, \eqref{0902:01}, and the choice of $r_0$, we get
\begin{align*}
 \int_{Q_+}w^{p_0}
 \geqslant\int_Ew^{p_0}
 \geqslant c_3r_0^N\left(\frac{w(z_0)}2\right)^{p_0}
\geqslant \frac{c_3}{4^{p_0}16^N}M^{p_0}
 =:c_4M^{p_0},
\end{align*}
which implies
\[
 {M\leqslant c_4^{-1/p_0}
       \left(\int_{Q_+}w^{p_0}\right)^{1/p_0}.}
\]
Finally, by the definition of $M$, we have
\[
 \sup_Qw
\leqslant\frac{M}{d_{\rm K}(Q,Q_+^c)^{N/p_0}}
 \leqslant \frac{c_4^{-1/p_0}}{d_{\rm K}(Q,Q_+^c)^{N/p_0}}
       \left(\int_{Q_+}w^{p_0}\right)^{1/p_0}.
\]
{Combining this estimate with the weak Harnack inequality
\eqref{eq:weak-harnack}, we complete the proof.}
\end{proof}

\subsection{H\"older regularity}
\label{sec:holder}

We now derive the interior H\"older estimate from the strong Harnack
inequality proved in Subsection \ref{sec:strong-harnack}. Recall the symmetric cylinder:
\[
 Q_0:=\{(t,x,v): |t|<1,\ |x|<1,\ |v|<1\}.
\]

\begin{lemma}\label{lem:Harnack_scal}
Under \eqref{AB2}, there are constants $A_0>4$, $t_0>0$, and $C_H>1$, depending only on $\Theta${,} such that for every $(z,r)\in\R^{1+2d}\times(0,1]$ and every nonnegative solution $w$
of $\cL_{a,b}w=0$ in
$\mathcal G_{z,r}\bigl(\delta_{A_0}(Q_0)\bigr)$,
\begin{equation}\label{eq:scaled-harnack}
{\sup_{\mathcal G_{z,r}(z_0\circ Q_0)}w
 \leqslant C_{\rm H}\inf_{\mathcal G_{z,r}(Q_0)}w,}
\end{equation}
where $z_0:=(t_0,0,0)\in \R^{1+2d}$.
\end{lemma}
\begin{proof}
{Fix}
\[
 z_{\rm L}:=(-9.5,0,0),\qquad z_{\rm R}:=(-0.5,0,0).
\]
{Then}
\[
 {\mathcal G_{z_{\rm L},1/2}(Q_0)\Subset\mathring Q_-,
 \qquad
 \mathcal G_{z_{\rm R},1/2}(Q_0)\Subset\mathring Q_+.}
\]
{Applying} Theorems \ref{thm:weak-harnack} and
\ref{thm:strong-harnack} with ${2\hbar}$ in place of $\hbar$, and {denoting} the
corresponding outer radius by $R_*>4${, we find} a constant $C_{\rm H}\geqslant2$, depending only on
$\Theta$, such that every nonnegative solution $u$ {of $\cL_{a,b}u=0$ in
$\delta_{R_*}(Q_0)$ satisfies}
\begin{equation}\label{0903:00}
 {\sup_{\mathcal G_{z_{\rm R},1/2}(Q_0)}u
 \leqslant C_{\rm H}\inf_{Q_-}u
 \leqslant C_{\rm H}
 \inf_{\mathcal G_{z_{\rm L},1/2}(Q_0)}u.}
\end{equation}
{Set}
\[
z_0:=\delta_2(z_{\rm L}^{-1}\circ z_{\rm R})
 =(36,0,0)=:(t_0,0,0),\quad
\widehat z:=(\delta_2z_{\rm L})^{-1}=(38,0,0).
\]
{By \eqref{KiDi2},} we have
\[
 {\mathcal G_{\widehat z,2}\circ\mathcal G_{z_{\rm L},1/2}
 =\mathcal G_{0,1},
 \quad
\mathcal G_{\widehat z,2}\circ\mathcal G_{z_{\rm R},1/2}
 =\mathcal G_{z_0,1}.}
\]
{Choose}
\[
 {A_0>\max\bigl\{4,\normK{\widehat z}+2R_*\bigr\}.}
\]
{Indeed, \eqref{Af3} and the homogeneity of $\normK{\cdot}$ give,
for every $z'\in\delta_{R_*}(Q_0)$,}
\[
 {\normK{\mathcal G_{\widehat z,2}(z')}
 \leqslant\normK{\widehat z}+2\normK{z'}
 <A_0.}
\]
{Consequently,}
\begin{equation}\label{eq:T-outer-containment}
 {\mathcal G_{\widehat z,2}\bigl(\delta_{R_*}(Q_0)\bigr)
 \subset\delta_{A_0}(Q_0).}
\end{equation}
Fix $(z,r)\in\R^{1+2d}\times(0,1]$ and let $w$ be as in the
statement.  Set
\[
 \mathcal G_{z,r}\circ\mathcal G_{\widehat z,2}
 =\mathcal G_{\widetilde z,\widetilde r},
 \qquad
 \widetilde z:=\mathcal G_{z,r}(\widehat z),
 \qquad
 \widetilde r:=2r
\]
By \eqref{eq:T-outer-containment} and Lemma \ref{lem_scaling_PDE},
$\widetilde w:=w\circ\mathcal G_{\widetilde z,\widetilde r}$ is a nonnegative solution in
$\delta_{R_*}(Q_0)$ with
\[
\widetilde a=a\circ\mathcal G_{\widetilde z,\widetilde r},
 \qquad
 \widetilde b=\widetilde r\,b\circ\mathcal G_{\widetilde z,\widetilde r}.
\]
Note that
\[
 {\mathcal G_{\widetilde z,\widetilde r}\circ\mathcal G_{z_{\rm L},2}
 =\mathcal G_{z,r},
 \quad
\mathcal G_{\widetilde z,\widetilde r}\bigl(\mathcal G_{z_{\rm R},2}(Q_0)\bigr)
 =\mathcal G_{z,r}(z_0\circ Q_0).}
\]
Applying \eqref{0903:00} to $\widetilde w$, we  prove
\eqref{eq:scaled-harnack}.
\end{proof}

{For every nonempty set $O\subset\R^{1+2d}$, define the oscillation
of $w$ on $O$ by}
\[
 {\osc_O w:=\sup_Ow-\inf_Ow.}
\]

{We next establish a local oscillation-decay estimate.}

\begin{lemma}[{Oscillation decay}]\label{lem:osc}
{Under \eqref{AB2}, there exist constants $\alpha\in(0,1)$ and
$C_{\rm osc}>0$, depending only on $\Theta$, with the
following property:  Let $\Omega\subset\R^{1+2d}$ be open, let $w$ be a
solution of $\cL_{a,b}w=0$ in $\Omega$, and suppose that for some $z\in\R^{1+2d}$,}
\[
\overline{\mathcal G_{z,1/2}(Q_0)}\Subset\Omega.
\]
{Then, for every $r\in(0,1/2]$,}
\begin{equation}\label{eq:two-scale-osc}
 {\osc_{\mathcal G_{z,r}(Q_0)}w
 \leqslant C_{\rm osc}r^\alpha
 \osc_{\mathcal G_{z,1/2}(Q_0)}w.}
\end{equation}
\end{lemma}

\begin{proof}
{Let $z_0,A_0$, and $C_{\rm H}$ be the constants from Lemma
\ref{lem:Harnack_scal}.  For $n=0,1,\ldots$, set}
\[
 {\rho_n:=A_0^{-n}/2,
 \quad
 E_n:=\mathcal G_{z,\rho_n}(z_0\circ Q_0),
 \qquad
 F_n:=\mathcal G_{z,\rho_n}(Q_0),}
\]
{and define}
\[
 {M_n:=\sup_{F_n}w,
 \qquad
 m_n:=\inf_{F_n}w.}
\]
{For every $n\geqslant1$, by \eqref{KiDi2}, we have}
\begin{equation}\label{eq:osc-outer-cylinder}
 {F_{n-1}
 =\mathcal G_{z,\rho_n}\bigl(\delta_{A_0}(Q_0)\bigr).}
\end{equation}
{Both $w-m_{n-1}$ and $M_{n-1}-w$ are nonnegative solutions in
$F_{n-1}$.  By Lemma \ref{lem:Harnack_scal}, we obtain}
\begin{align*}
 {\sup_{E_n}(w-m_{n-1})}
 &{\leqslant C_{\rm H}(m_n-m_{n-1}),}\\
 {\sup_{E_n}(M_{n-1}-w)}
 &{\leqslant C_{\rm H}(M_{n-1}-M_n).}
\end{align*}
{Adding the two inequalities, we get}
\[
 {\osc_{F_{n-1}}w+\osc_{E_n}w
 \leqslant C_{\rm H}
 \bigl(\osc_{F_{n-1}}w-\osc_{F_n}w\bigr).}
\]
{Since $\osc_{E_n}w\geqslant0$, it follows that}
\[
 {\osc_{F_n}w
 \leqslant\left(1-C^{-1}_{\rm H}\right)\osc_{F_{n-1}}w.}
\]
{Set}
\[
\theta:=1-C^{-1}_{\rm H}\in[1/2,1),
 \quad
 \alpha:=-\log\theta/\log A_0.
\]
{Since $A_0>4$ and $\theta\geqslant1/2$, we have
$\alpha\in(0,1)$.  By iteration, we get}
\[
 {\osc_{F_n}w\leqslant\theta^n\osc_{F_0}w
 \quad\text{for every }n\in\mathbb N_0.}
\]
{For $r\in(0,1/2]$, let}
\[
 {n_r:=\left\lfloor
 \frac{\log(1/(2r))}{\log A_0}\right\rfloor\in\mathbb N_0.}
\]
{Then}
\[
 {\frac{A_0^{-(n_r+1)}}2<r
 \leqslant\frac{A_0^{-n_r}}2<A_0r.}
\]
{Consequently, $\mathcal G_{z,r}(Q_0)\subset F_{n_r}$ and}
\begin{align*}
 {\osc_{\mathcal G_{z,r}(Q_0)}w}
 &{\leqslant\osc_{F_{n_r}}w
 \leqslant\theta^{n_r}\osc_{F_0}w}\\
 &{=(A_0^{-n_r})^\alpha\osc_{F_0}w
 \leqslant(2A_0r)^\alpha
 \osc_{\mathcal G_{z,1/2}(Q_0)}w.}
\end{align*}
{This proves \eqref{eq:two-scale-osc} with
$C_{\rm osc}:=(2A_0)^\alpha$.}
\end{proof}

\begin{theorem}[Interior kinetic H\"older estimate]
\label{thm:holder}
{Under \eqref{AB2}, there exist constants $\alpha\in(0,1)$ and
$C>0$, depending only on $\Theta$, such that
for any solution $w$ of $\cL_{a,b}w=0$ in $Q_0$,}
\begin{equation}\label{eq:interior-holder}
 {\sup_{z_1\ne z_2\in\mathcal G_{0,1/2}(Q_0)}
 \frac{|w(z_1)-w(z_2)|}{d_{\rm K}(z_1,z_2)^\alpha}
 \leqslant C\osc_{Q_0}w.}
\end{equation}
\end{theorem}

\begin{proof}
{Fix distinct points
$z_1,z_2\in\mathcal G_{0,1/2}(Q_0)$, and write
$z_1=(t_1,x_1,v_1)$.  The explicit formula for $\mathcal G$ shows that}
\[
 {\overline{\mathcal G_{z_1,1/2}(Q_0)}\Subset Q_0.}
\]
{Indeed, points in this closed cylinder satisfy}
\[
 {|t|\leqslant|t_1|+\tfrac14<\tfrac12,
 \qquad
 |x|\leqslant|x_1|+\tfrac18+\tfrac14|v_1|<\tfrac38,
 \qquad
 |v|\leqslant|v_1|+\tfrac12<1.}
\]
{Write $d_{12}:=d_{\rm K}(z_1,z_2)$.  If $d_{12}\geqslant1/4$,
then}
\[
 {|w(z_1)-w(z_2)|
 \leqslant\osc_{Q_0}w
 \leqslant4^\alpha d_{12}^\alpha\osc_{Q_0}w.}
\]
{Suppose instead that $d_{12}\in(0,1/4)$, and set
$\rho:=2d_{12}\in(0,1/2)$.  By the definition and homogeneity of the
kinetic distance,}
\[
 {\normK{\delta_{\rho^{-1}}(z_1^{-1}\circ z_2)}
 \leqslant\tfrac{d_{12}}\rho=\tfrac12<1.}
\]
Hence
$\delta_{\rho^{-1}}(z_1^{-1}\circ z_2)\in Q_0$  and
\[
 {z_1,z_2\in\mathcal G_{z_1,\rho}(Q_0).}
\]
By Lemma \ref{lem:osc} with $\Omega=Q_0$, we get
\begin{align*}
 {|w(z_1)-w(z_2)|}
 &{\leqslant\osc_{\mathcal G_{z_1,\rho}(Q_0)}w\leqslant C_{\rm osc}\rho^\alpha
 \osc_{\mathcal G_{z_1,1/2}(Q_0)}w
 \leqslant2^\alpha C_{\rm osc}d_{12}^\alpha\osc_{Q_0}w.}
\end{align*}
{Combining the two cases proves \eqref{eq:interior-holder}.}
\end{proof}

\section{Proof of Theorem \ref{main}}\label{sec:proof-main}

We now deduce the main theorem from the zero-order-free estimates proved
above.  The weak Harnack inequality follows from a time-exponential
transformation, whereas the strong Harnack and H\"older estimates follow
from a one-dimensional lifting argument.

\begin{proof}[Proof of Theorem~\ref{main}]
(i) Suppose that $w$ is a nonnegative strong supersolution of
$\cL_{a,b,c}$ in $\Omega$, and define
\[
 u(t,x,v):=\e^{-\hbar t}w(t,x,v).
\]
Since $c+\hbar\geqslant0$, we have
\begin{equation}\label{eq:time-exponential-reduction}
 \cL_{a,b}u
 =\e^{-\hbar t}
 \bigl(\cL_{a,b,c}w-(c+\hbar)w\bigr)
 \leqslant0
 \qquad\text{a.e. in }\Omega.
\end{equation}
Thus $u$ is a nonnegative strong supersolution of $\cL_{a,b}$.  By Theorem~\ref{thm:weak-harnack}, we have
\[
 \left(\fint_{Q_+}u^{p_0}\right)^{1/p_0}
 \leqslant C_0\inf_{Q_-}u,
\]
which in turn implies
\[
 \left(\fint_{Q_+}w^{p_0}\right)^{1/p_0}
 \leqslant \e^{10\hbar}C_0\inf_{Q_-}w.
\]
This proves assertion~{\rm(i)}.

{\rm(ii)} and~{\rm(iii)}.  
Introduce an
additional position--velocity pair $(\xi,\eta)\in\R\times\R$ and write
\[
 \widehat x:=(x,\xi),\qquad
 \widehat v:=(v,\eta).
\]
On the $(d+1)$-dimensional kinetic space, define
\[
 \widehat w(t,\widehat x,\widehat v)
 :=\e^\eta w(t,x,v),
 \qquad
 \widehat a:=
 \begin{pmatrix}
  a&0\\
  0&1
 \end{pmatrix},
 \qquad
 \widehat b:=(b,c-1),
\]
where the coefficients on the right-hand side are evaluated at
$(t,x,v)$.  Let
\[
 \widehat\Omega
 :=\{(t,(x,\xi),(v,\eta)):(t,x,v)\in\Omega\}.
\]
Since $\widehat w$ is independent of $\xi$, by direct calculation, we have
\begin{equation}\label{eq:dimension-lift}
 \cL_{\widehat a,\widehat b}\widehat w
 =\e^\eta\cL_{a,b,c}w=0.
\end{equation}
Thus $\widehat w$ is a strong solution of $\cL_{\widehat a,\widehat b}$ in $\widehat\Omega$.
Moreover,
\[
 \widehat\lambda I_{d+1}
 \preceq\widehat a
 \preceq\widehat\Lambda I_{d+1},
 \qquad
 |\widehat b|
 \leqslant\sqrt{\hbar^2+(\hbar+1)^2}
 =:\widehat\hbar,
\]
where
\[
 \widehat\lambda:=\min\{\lambda,1\},\quad
 \widehat\Lambda:=\max\{\Lambda,1\}.
\]
Hence the lifted equation satisfies the zero-order-free structural
conditions with parameter tuple 
\[
 \widehat\Theta
 :=(d+1,\widehat\lambda,\widehat\Lambda,\widehat\hbar).
\]
Denote by $\widehat Q_0$ and $\widehat Q_\pm$ the reference cylinders in
dimension $d+1$, and let $\widehat d_{\rm K}$ be the corresponding kinetic
distance.  
Let $\Omega=\delta_{R_m}(Q_0)$ and $\widehat\Omega=\delta_{\widehat R_m}(\widehat Q_0)$, where $R_m$ 
and $\widehat R_m$ are the radius given by
Theorem~\ref{thm:weak-harnack} for the parameter tuples $\Theta$ and $\widehat\Theta$. 
Thus, 
$$
Q_0\cup Q_+\cup Q_-\Subset\Omega,\quad \widehat Q_0\cup \widehat Q_+\cup \widehat Q_-\Subset\widehat\Omega.
$$ 
The projection of
$\delta_{\widehat R}(\widehat Q_0)$ onto the $(t,x,v)$ variables is
contained in $\delta_{\widehat R}(Q_0)\subset\Omega$; hence
\[
 \delta_{\widehat R}(\widehat Q_0)\subset\widehat\Omega.
\]
Let
\[
 \iota(t,x,v):=(t,(x,0),(v,0))
\]
be the zero-section embedding.  For
$\delta:=d_{\rm K}(Q,Q_+^c)>0$, by
\eqref{eq:group}--\eqref{eq:metric}, we have for some $\kappa_0\in(0,\frac14)$,
\begin{equation}\label{eq:lifted-separation}
 \widehat d_{\rm K}\bigl(\iota(Q),\widehat Q_+^c\bigr)
 \geqslant \kappa_0\min\{\delta,\delta^{3/2},1\}{>0.}
\end{equation}
Indeed, fix $z=(t,x,v)\in Q$. By the definition of $\delta$, we have
\[
 t+1\geqslant\delta^2,\qquad -t\geqslant\delta^2,
 \qquad 1-|v|\geqslant\delta,
 \qquad 1-|x|\geqslant\delta^3.
\]
Set $\rho:=\kappa_0\min\{\delta,\delta^{3/2},1\}$ and suppose that
$\widehat z'=(s,(y,\xi),(\nu,\eta))$ satisfies
$\widehat d_{\rm K}(\iota(z),\widehat z')<\rho$.  Then
\[
 |s-t|<\rho^2,\qquad
 |(\nu-v,\eta)|<\rho,\qquad
 |(y-x-(s-t)v,\xi)|<\rho^3.
\]
Consequently, for $\kappa_0$ {chosen} sufficiently small,
\[
 -1<s<0,\quad
 |(\nu,\eta)|<1,\quad
 |(y,\xi)|
 \leqslant |x|+\rho^3+\rho^2|v|<1.
\]
Thus the kinetic
ball of radius $\rho$ about $\iota(z)$ is contained in
$\mathring{\widehat Q}_+$, which proves \eqref{eq:lifted-separation}.

If $w$ is a nonnegative strong solution, then $\widehat w\geqslant0$.
By Theorem~\ref{thm:strong-harnack} in dimension $d+1$ with $\iota(Q)\Subset\mathring{\widehat Q}_+$,
 we get
\begin{align*}
 \sup_Qw
 &=\sup_{\iota(Q)}\widehat w
 \leqslant
 C\bigl(\widehat\Theta,
   \widehat d_{\rm K}(\iota(Q),\widehat Q_+^c)\bigr)
 \inf_{\widehat Q_-}\widehat w \leqslant C_1\inf_{Q_-}w.
\end{align*}
Here we used $\iota(Q_-)\subset\widehat Q_-$ and
$\widehat w\circ\iota=w$.  This proves assertion~{\rm(ii)}.

Finally, let $w$ be any strong solution. Noting that
$\iota(\delta_{1/2}(Q_0))\subset\delta_{1/2}(\widehat Q_0)$,
by applying Theorem~\ref{thm:holder} to $\widehat w$ in $\widehat Q_0$, {we obtain,
for the exponent $\alpha=\alpha(\widehat\Theta)\in(0,1)$,}
\[
 \sup_{z_1,z_2\in\delta_{1/2}(Q_0), z_1\ne z_2}
 \frac{|w(z_1)-w(z_2)|}{d_{\rm K}(z_1,z_2)^\alpha}
 \leqslant
 \widehat C\,\osc_{\widehat Q_0}\widehat w.
\]
Since the projection of $\widehat Q_0$ is contained in $Q_0$ and
$|\eta|<1$ there,
\[
 \osc_{\widehat Q_0}\widehat w
 \leqslant2\|\widehat w\|_{L^\infty(\widehat Q_0)}
 \leqslant2\e\|w\|_{L^\infty(Q_0)}.
\]
{This proves assertion~{\rm(iii)}.} 
\end{proof}

\appendix
\section{Auxiliary results}

\subsection{Stochastic estimates}
We introduce the following lower bound for the probability that a nondegenerate martingale
remains in a ball.
\begin{proposition}[Small-ball survival]\label{prop:sb}
Let $(M_t)_{t{\geqslant}0}$ be a continuous $\R^d$-valued martingale on a filtered probability
space, with $M_0=0$.  Assume that its quadratic variation is absolutely continuous and
\begin{equation}\label{eq:martingale-ellipticity}
2\lambda I_d\preceq
\frac{\dd\qv{M}_t}{\dd t}
\preceq2\Lambda I_d
\qquad \dd t\otimes\dd\Pp\text{-a.e.},
\end{equation}
where $0<\lambda\leqslant\Lambda<\infty$.  Then 
there is a constant $\mu_1=\mu_1(\lambda,\Lambda,d)>0$ such that for every $T,R>0$,
\begin{equation}\label{eq:sb}
\Pp\Bigl(\sup_{0\leqslant s\leqslant T}|M_s|<R\Bigr)
{\geqslant}
\exp\!\left(-\mu_1\frac{T}{R^2}\right).
\end{equation}
\end{proposition}

\begin{proof}
Define the Pucci operator on $\Ssym^d$ by
\begin{equation}\label{eq:infrep}
\Pucci(A)
=
{\inf\bigl\{\tr(QA):Q\in\Ssym^d,\
\lambda I_d\preceq Q\preceq\Lambda I_d\bigr\}.}
\end{equation}
Consider the eigenvalue problem
\begin{equation}\label{eq:eigen}
\Pucci(\nabla^2\varphi)=-\mu_1\varphi\quad\text{in }B_1,
\qquad
\varphi=0\quad\text{on }\partial B_1.
\end{equation}
By the standard principal-eigenvalue theory for uniformly elliptic, positively homogeneous
operators (\cite{BEQ05,QS08,EFQ10}), there is a principal eigenvalue $\mu_1\in(0,\infty)$ and a positive eigenfunction
\[
\varphi\in C^{2,\alpha}_{\mathrm{loc}}(B_1)\cap C(\overline B_1)
\]
for some $\alpha\in(0,1)$.  The positive eigenfunction is unique up
to a positive factor. By the rotational invariance, 
$\varphi(x)=\phi(|x|)$ is radial. Furthermore, by the strong maximum
principle, for every $r\in(0,1)$, the minimum of $\varphi$ on
$\overline{B_r}$ is attained only on $\partial B_r$. Hence, $\phi$ is strictly decreasing.
Thus we may normalize
it by
\begin{equation}\label{eq:eigenfunction-normalization}
\varphi(0)=\|\varphi\|_{L^\infty(B_1)}=1.
\end{equation}

{For $R>0$, set $\widetilde M_t:=R^{-1}M_{R^2t}$.}
Then $\qv{\widetilde M}_t=R^{-2}\qv{M}_{R^2t}$ and
\[
\Pp\Bigl(\sup_{0\leqslant s\leqslant T}|M_s|<R\Bigr)
=
\Pp\Bigl(\sup_{0\leqslant u\leqslant T/R^2}|\widetilde M_u|<1\Bigr).
\]
It is therefore enough to prove the result for $R=1$.
Let
\[
{\tau:=\inf\{t\geqslant0:|M_t|\geqslant1\},}
\qquad
{\tau_n:=\inf\{t\geqslant0:|M_t|\geqslant1-1/n\},\quad
n\geqslant2.}
\]
By continuity, we have $\tau_n\uparrow\tau$.  For fixed $T>0$, define
\[
Y_t^{(n)}
:=
{\e^{\mu_1(t\wedge\tau_n)}}
\varphi(M_{t\wedge\tau_n}),
\qquad 0\leqslant t\leqslant T.
\]
Since $M_{t\wedge\tau_n}\in\overline B_{1-1/n}$ and
$\varphi\in C^2(B_1)$, by It\^o's formula, we have up to $\tau_n$,
\begin{align*}
\dd\bigl({\e^{\mu_1t}}\varphi(M_t)\bigr)
={}&{\e^{\mu_1t}}
\left[\mu_1\varphi(M_t)\dd t
+\tfrac12\nabla^2\varphi(M_t):\dd \qv{M}_t\right]
+{\e^{\mu_1t}}\nabla\varphi(M_t)\cdot\dd M_t.
\end{align*}
By \eqref{eq:eigen}, we have
\[
\tfrac12\nabla^2\varphi(M_t):\tfrac{\dd \qv{M}_t}{\dd t}
{\geqslant}
\Pucci(\nabla^2\varphi(M_t))
=-\mu_1\varphi(M_t).
\]
In particular,
$$
\E\left[
{\e^{\mu_1(T\wedge\tau_n)}}
\varphi(M_{T\wedge\tau_n})
\right]
\geqslant\varphi(0)=1.
$$
Since $\tau_n\uparrow\tau$ and the stopped integrands are uniformly
bounded, by the dominated convergence we get
\[
{\begin{aligned}
1\leqslant
\E\left[\e^{\mu_1(T\wedge\tau)}
\varphi(M_{T\wedge\tau})\right]=\e^{\mu_1T}\E\bigl[\varphi(M_T)\1_{\{\tau>T\}}\bigr]
\leqslant
\e^{\mu_1T}\Pp(\tau>T).
\end{aligned}}
\]
{Here we used $\varphi=0$ on $\partial B_1$ and
$0\leqslant\varphi\leqslant1$.  Since}
\[
 {\{\tau>T\}
 =\left\{\sup_{0\leqslant s\leqslant T}|M_s|<1\right\},}
\]
{we obtain \eqref{eq:sb}.}
\end{proof}

\begin{lemma}
\label{lem:tube}
{Let $h>0$ and let $\phi:[0,h]\to\R^d$ be absolutely continuous with}
$\phi(0)=0$, and
\[
 \mathcal A_\phi:=\int_0^h|\dot\phi_s|^2\dd s<\infty.
\]
{Let $\hbar\geqslant0$.  Then, for every progressively measurable}
$(\mathbf A,\boldsymbol\beta)$ satisfying
$\lambda I_d\preceq\mathbf A_s\preceq\Lambda I_d$ and
$|\boldsymbol\beta_s|\leqslant \hbar$, and every $R>0$,
\begin{align}
 \Pp\bigg(
 {\sup_{0\leqslant s\leqslant h}}\left|
 \int_0^s\boldsymbol\beta_r\dd r
 +\int_0^s\sqrt{2\mathbf A_r}\dd W_r-\phi_s\right|<R
 \bigg)
 \geqslant
 \exp\left(-\left(
 \frac{\sqrt{\mu_1h}}R
 +\sqrt\frac{\mathcal A_\phi+\hbar^2h}{2\lambda}
 \right)^2\right),
 \label{eq:tube}
\end{align}
where $\mu_1$ is the same as in Proposition \ref{prop:sb}.
\end{lemma}

\begin{proof}
Write $\boldsymbol\sigma_s:=\sqrt{2\mathbf A_s}$.  Then
\[
 \psi_s:=\phi_s-\int_0^s\boldsymbol\beta_r\dd r
\]
{is an adapted absolutely continuous process.  Since
$\|\boldsymbol\sigma_s^{-1}\|\leqslant(2\lambda)^{-1/2}$ and
$|u-v|^2\leqslant2|u|^2+2|v|^2$,}
\begin{equation}
{\int_0^h|\boldsymbol\sigma_s^{-1}\dot\psi_s|^2\dd s
 \leqslant \frac{\mathcal A_\phi+\hbar^2h}{\lambda}.}
 \label{eq:girsanov-energy}
\end{equation}
By Novikov's criterion, we may define a probability measure $\mathbb Q^\phi$ by
\[
 {\frac{\dd\mathbb Q^\phi}{\dd\Pp}}
 =\exp\left({\int_0^h\boldsymbol\sigma_s^{-1}\dot\psi_s\cdot\dd W_s}
 -\frac12\int_0^h|\boldsymbol\sigma_s^{-1}\dot\psi_s|^2\dd s\right).
\]
By Girsanov's theorem,
\[
 W_s^\phi
 :=
 W_s-\int_0^s\boldsymbol\sigma_r^{-1}\dot\psi_r\dd r
\]
{is a Brownian motion under $\mathbb Q^\phi$, which is equivalent to
$\Pp$.  Thus, under $\mathbb Q^\phi$,}
\[
 M^\phi_s:=
 \int_0^s\boldsymbol\beta_r\dd r
 +\int_0^s\boldsymbol\sigma_r\dd W_r-\phi_s=
 \int_0^s\boldsymbol\sigma_r\dd W_r-\psi_s
 =
\int_0^s\boldsymbol\sigma_r\dd W_r^\phi
\]
is a continuous martingale starting from zero and satisfying
\[
 2\lambda I_d
 \preceq
 \frac{\dd\qv{M^\phi}_s}{\dd s}
 \preceq
 2\Lambda I_d.
\]
By Proposition \ref{prop:sb}, we have
\begin{align}\label{Aw1}
 \mathbb Q^\phi\bigg(
 {\sup_{0\leqslant s\leqslant h}}|M^\phi_s|<R
 \bigg)\geqslant {\e^{-\mu_1h/R^2}}.
\end{align}
For $p>1$, the exponential-martingale identity and
\eqref{eq:girsanov-energy} yield
\[
 \E_{\Pp}\left(\frac{\dd\mathbb Q^\phi}{\dd\Pp}\right)^p
 \leqslant\exp\left(\frac{p(p-1)(\mathcal A_\phi+\hbar^2h)}{2\lambda}\right).
\]
Now, by H\"older's inequality, we get
\[
\mathbb Q^\phi\bigg(
 {\sup_{0\leqslant s\leqslant h}}|M^\phi_s|<R
 \bigg)\leqslant \Pp\bigg(
 {\sup_{0\leqslant s\leqslant h}}|M^\phi_s|<R
 \bigg)^{1-1/p}
 \exp\left(\frac{(p-1)(\mathcal A_\phi+\hbar^2h)}{2\lambda}\right),
\]
which together with \eqref{Aw1} yields
\[
 \Pp\bigg(
 {\sup_{0\leqslant s\leqslant h}}|M^\phi_s|<R
 \bigg)
 \geqslant\exp\left(-\frac{p\mu_1h}{(p-1)R^2}
 -\frac{p(\mathcal A_\phi+\hbar^2h)}{2\lambda}\right).
\]
With $a:=\mu_1h/R^2$ and
$c:=(\mathcal A_\phi+\hbar^2h)/(2\lambda)$,
{if $c>0$, the exponent is minimized by
$p-1=\sqrt{a/c}$.  If $c=0$, we let $p\to\infty$.  In either case, its
infimum over $p>1$ is $(\sqrt a+\sqrt c)^2$.}
This proves \eqref{eq:tube}.
\end{proof}

\subsection{A kinetic ink-spots lemma}
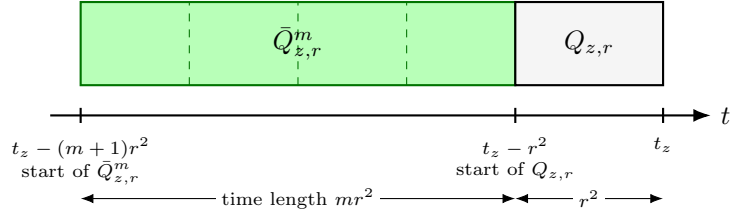
\begin{figure}[htbp]

\centering

\begin{tikzpicture}[x=1.15cm,y=1cm,>=Latex]
\fill[green!28] (0,0.65) rectangle (5,1.75);
\draw[green!45!black,thick] (0,0.65) rectangle (5,1.75);

\draw[dashed,green!45!black] (1.25,0.65) -- (1.25,1.75);
\draw[dashed,green!45!black] (2.50,0.65) -- (2.50,1.75);
\draw[dashed,green!45!black] (3.75,0.65) -- (3.75,1.75);

\fill[gray!8] (5,0.65) rectangle (6.7,1.75);
\draw[thick] (5,0.65) rectangle (6.7,1.75);

\node at (2.5,1.2) {${\bar{Q}_{z,r}}^m$};
\node at (5.85,1.2) {$Q_{z,r}$};

\draw[->,thick] (-0.35,0.25) -- (7.25,0.25)
  node[right] {$t$};

\foreach \x in {0,5,6.7}
  \draw[thick] (\x,0.15) -- (\x,0.35);

\node[below,align=center,font=\scriptsize] at (0,0.08)
  {$t_z-(m+1)r^2$\\
   start of ${\bar{Q}_{z,r}}^m$};

\node[below,align=center,font=\scriptsize] at (5,0.08)
  {$t_z-r^2$\\
   start of $Q_{z,r}$};

\node[below,font=\scriptsize] at (6.7,0.08)
  {$t_z$};

\draw[<->] (0,-0.85) -- (5,-0.85)
  node[midway,fill=white,inner sep=2pt,font=\scriptsize]
  {time length $mr^2$};

\draw[<->] (5,-0.85) -- (6.7,-0.85)
  node[midway,fill=white,inner sep=2pt,font=\scriptsize]
  {$r^2$};

\end{tikzpicture}

\caption{Time arrangement of the predecessor stack
${\bar{Q}_{z,r}}^m$ and the cylinder $Q_{z,r}$.}
\label{fig:predecessor-stack}
\end{figure}

Recall that
$$
Q_+:=\{-1<t\leqslant 0,\, |x|\leqslant1,\,  |v|\leqslant 1\},
$$
and for any $z=(t_z,x_z,v_z)\in \mathbb{R}^{1+2d}$ and $r>0$,
$$
Q_{z,r}:=\mathcal G_{z,r}(Q_+)=\{-r^2<t-t_z\leqslant0,\, |x-x_z-r^2v_z|\leqslant r^3,\, |v-v_z|\leqslant r \}.
$$
For any $m\geqslant1$, we define (see Figure \ref{fig:predecessor-stack})
$$
{\bar{Q}_{z,r}}^m:=\{-(m+1)r^2\leqslant t-t_z<-r^2,\, |x-x_z-(t-t_z)v_z|\leqslant (m+2)r^3,\, |v-v_z|\leqslant r \}.
$$

\begin{theorem}[Ink-spots theorem with leakage] \label{thm:ink-spots}
Let $E, F$ be two measurable {subsets} of $\R^{1+2d}$ with
$$
E \subset F\cap Q_+. 
$$
Let $\mu\in(0,1)$ be a given constant. Define a set
$$
\mathcal E_{\mu}:=\Big\{(z,r)\in\R^{1+2d}\times(0,1): Q_{z,r} \subset Q_+,\ |Q_{z,r}\setminus E| \leqslant \mu |Q_{z,r}| \Big\}.
$$
Suppose that there are $r_0\in(0,1)$ and $m\in\mathbb N$ such that for any $(z,r)\in \mathcal E_{\mu}$,
\begin{align}\label{Af9}
{\bar{Q}_{z,r}}^m \subset F\mbox{ and }r<r_0.
\end{align}
Then there {exist} two constants $c_0\in(0,1)$ and $C_0\geqslant 1${,} depending only on $d${,} such that 
\begin{equation}
|E| \leqslant \frac{m+1}{m}(1-c_0\mu)\left(|F \cap Q_+| + C_0mr_0^{2}\right).
\end{equation}
\end{theorem}
\begin{proof}
In \cite[Corollary~10.2]{ImbertSilvestre20}, the stacked cylinder is obtained by translating the original cylinder $m$ steps backward in time, whereas our stacked cylinder is translated $m$ steps forward in time. The time reversal $(t,x,v)\mapsto(-t,x,-v)$ interchanges these two configurations while preserving all geometric and measure-theoretic properties. The result therefore follows directly from \cite[Corollary~10.2]{ImbertSilvestre20}.
\end{proof}

{\bf Acknowledgements.}
This work is supported by the NNSF of China (No. 12595282, 12601259) and the National Key R\&D
Program of China (No. 2023YFA1010103). 
We also acknowledge the use of AI tools (ChatGPT 5.6 Pro) for assistance in collecting relevant literature and for inspiring discussions, in particular for suggesting the use of  $(\mathsf Y_t, \mathsf V_t)$  in proving Krylov's estimate. In addition, AI tools were used to improve the language and to assist in preparing the two figures. All proofs in this paper were developed and written entirely by the authors.

\end{document}